\documentclass[a4paper,10pt,reqno]{amsart}

\usepackage{lineno}

\usepackage{BOONDOX-cal}
\usepackage{amssymb}
\usepackage{amsmath}
\usepackage[mathscr]{eucal}
\usepackage{mathrsfs}
\usepackage{xspace} 
\usepackage{tabu}
\usepackage{thmtools}
\usepackage[latin1]{inputenc}
\usepackage[english]{babel}
\usepackage[all,2cell]{xy}
\usepackage{xcolor}
\usepackage{graphicx}
\usepackage{wasysym} 
\usepackage{hyperref} 

\usepackage{pdflscape}
\usepackage{caption}
\usepackage{multicol}
\usepackage{longtable}

\makeatletter
\newcommand{\leqnomode}{\tagsleft@true}
\newcommand{\reqnomode}{\tagsleft@false}
\makeatother

\usepackage{array}
\newcolumntype{L}[1]{>{\raggedright\let\newline\\\arraybackslash\hspace{0pt}}m{#1}}
\newcolumntype{C}[1]{>{\centering\let\newline\\\arraybackslash\hspace{0pt}}m{#1}}
\newcolumntype{R}[1]{>{\raggedleft\let\newline\\\arraybackslash\hspace{0pt}}m{#1}}

\definecolor{sg}{HTML}{df78ef}
\definecolor{sg1}{HTML}{ab47bc}
\definecolor{sg2}{HTML}{790e8b}

\usepackage{tikz}
\usetikzlibrary{mindmap,backgrounds, calc}
\usetikzlibrary{matrix,chains,scopes,positioning,arrows,fit}
\usetikzlibrary{positioning, shapes, patterns}
\usetikzlibrary{automata} 
\usetikzlibrary{shapes.geometric, arrows,arrows.meta}
\usetikzlibrary{positioning,decorations.markings}
\usetikzlibrary{cd}
\usepackage{hyperref}
\usepackage[acronym]{glossaries}

\usepackage{imakeidx}
\makeindex

\makeatletter
\newcommand*\bigcdot{\mathpalette\bigcdot@{.5}}
\newcommand*\bigcdot@[2]{\mathbin{\vcenter{\hbox{\scalebox{#2}{$\m@th#1\bullet$}}}}}
\makeatother

\usepackage{stmaryrd} 
\UseAllTwocells

\newcommand{\at}[1]%
            {\ensuremath{\protect\underline{\mathbf{#1}}}} 

\newcommand{\op}[1]{\ensuremath{\operatorname{#1}}}        

\newcommand{\h}[1][]                                       
 {\ifthenelse{\boolean{mmode}}%
  {$\mathrm{h}$}%
  {h\nobreakdash#1\hspace{0pt}}}

\newcommand{\comp}{\circ}          
\newcommand{\adcomp}%
  {\overset{\operatorname{ad}}{\comp}} 
\newcommand{\funcomp}%
  {\overset{\operatorname{fn}}{\comp}}

\newcommand{\sccat}
{\mathbin{\kern-1pt\raisebox{6pt}{.}\kern-5pt
\downarrow\kern-5pt\raisebox{6pt}{.}\kern-1pt}}

\newcommand{\parrow}[1]
   {\underset{{\displaystyle \raisebox{5pt}%
   {$\longleftarrow$}}}{\op{#1}}{\,}}
\newcommand{\iarrow}[1]
   {\underset{{\displaystyle \raisebox{5pt}%
   {$\longrightarrow$}}}{\op{#1}}{\,}}

\newcommand{\rest}%
{\mathnormal{\restriction}}        

\newcommand{\brel}{\ensuremath{\xymatrix{{}\arity@{{*}{-}{*}}[r] & {}}}}

\newcommand{\nseq}[3]{\xymatrix@1@C=16pt{#1 \arity@{>}[r]_-{\scriptscriptstyle{#2}} & #3 }}

\NoCompileMatrices              
\xymatrixrowsep = {8ex}         
\xymatrixcolsep = {10ex}        

\newdir{<= }{:a(+25)\dir{-}*:a(-25)\dir{-}*!/:a(85)+1.8pt/:a(-25)\dir{-}}
\newdir{ = }{*@{=}}
\newdir{+>}{@{|}*@{-}!/-6.5pt/@{>}*!/-13pt/{}}
\newdir{ +}{{}*!/-5pt/@{+}}
\newdir{ >}{{}*!/-10pt/@{>}}
\newdir{.-}{{}*:a(-90)h!/2.1pt/{.}}
\newdir{.->}{:a(-90)h!/2.5pt/{.}*!/4pt/\dir{-}*!/0pt/\dir{-}*%
!/-2pt/\dir{-}*!/-8pt/\dir{>}}
\newdir{-.->}%
{:a(-90)h!/2pt/{.}*!/8pt/\dir{-}*!/4pt/\dir{-}*!/0pt/\dir{-}*%
!/-5pt/\dir{-}*!/-11pt/\dir{>}}
\newdir{~>}{!/4.6pt/\dir{~}*!/1pt/\dir{-}*!/-5pt/\dir{>}}
\newdir{~~>}{!/7pt/\dir{~}*!/0pt/\dir{~}*!/-5pt/\dir{>}}
\newdir{=~>}{*!/1pt/@2{~}*!/6pt/{}*!/-6pt/@2{>}}
\newdir{<~~>}{!/21pt/\dir{<}*!/21pt/\dir{-}!/12pt/\dir{~}!/5pt/\dir{~}*!/1pt/\dir{-}*!/-3pt/\dir{>}}
\newdir{=>}{!/5pt/\dir{=}!/2.5pt/\dir{=}*!/-5pt/\dir2{>}*!/7pt/\dir{ }}
\newdir{==>}{!/-2pt/\dir{=}!/-6pt/\dir{=}%
 !/-10pt/\dir{=}*!/-18pt/\dir2{>}}
\newdir{<==>}{!\dir2{<}!/-2pt/\dir{=}!/-6pt/\dir{=}%
 !/-10pt/\dir{=}*!/-17pt/\dir2{>}}
\newdir{< }{{}*!/12pt/\dir{<}}
\newdir{-o}{!/1.3pt/{}*!/0pt/{}@{o}*!/-5pt/{}}

\newsavebox{\xymor}  
\newsavebox{\xymon}  
\newsavebox{\xyepi}  
\newsavebox{\xytn}   
\newsavebox{\xyrel}  
\newsavebox{\xycel}  
\newsavebox{\xymdf}  
\newsavebox{\xyumor} 
\newsavebox{\xydmor} 
\newsavebox{\xyomor} 
\newsavebox{\xyemor} 

\newcommand{\xynode}{\makebox[0ex]{}}
\savebox{\xymor}{\ensuremath{%
\xymatrix@1@C=19pt{\xynode \ar@{>}[r] & \xynode }}}
\savebox{\xymon}{\ensuremath{%
\xymatrix@1@C=19pt{\xynode \ar@{{ +}{-}{>}}[r] & \xynode }}}
\savebox{\xyepi}{\ensuremath{%
\xymatrix@1@C=19pt{\xynode \ar@{{}{-}{+>}}[r] & \xynode }}}
\savebox{\xytn}{\ensuremath{%
\xymatrix@1@C=19pt{\xynode \ar[r]|(.44){\object@{.-}} & \xynode
}}}
\savebox{\xyrel}{\ensuremath{%
\xymatrix@1@C=19pt{\xynode \ar@{{}{-}{-o}}[r] & \xynode }}}
\savebox{\xycel}{\ensuremath{%
\xymatrix@1@C=19pt{\xynode \ar@{=>}[r] & \xynode }}}
\savebox{\xymdf}{\ensuremath{%
\xymatrix@1@C=16pt{\xynode \ar@{}[r]|{\dir{~>}} & \xynode}}}
\savebox{\xyumor}{\ensuremath{%
\xymatrix@1@C=19pt{\xynode \ar@{{}{-}^{>}}[r] & \xynode }}}
\savebox{\xydmor}{\ensuremath{%
\xymatrix@1@C=19pt{\xynode \ar@{{}{-}_{>}}[r] & \xynode }}}
\savebox{\xyomor}{\ensuremath{%
\xymatrix@1@C=19pt{\xynode \ar@{{}{-}^{< }}[r] & \xynode }}}
\savebox{\xyemor}{\ensuremath{%
\xymatrix@1@C=19pt{\xynode \ar@{{ >}{-}{>}}[r] & \xynode }}}

\newcommand{\mor}{\usebox{\xymor}}    

\newcommand{\functor}[9]{
 \xymatrix{
    #4 \save[]+<0ex,5ex>*+{#1}="1"  \restore
      \arity[d]_{#6}  \arity@{}[rd]|{\longmapsto}
  & #5 \save[]+<0ex,5ex>*+{#3}="3"  \restore
      \arity[d]^{#7}
  \\
   #8 & #9 \arity "1";"3"^-{#2} } }
\newcommand{\functornd}[9]{
 \xymatrix{
    #4 \save[]+<0ex,5ex>*+{#1}="1"  \restore
      \arity[d]_{#6}  \arity@{}[rd]|{\longmapsto}
  & #5 \save[]+<0ex,5ex>*+{#3}="3"  \restore
  \\
   #8 & #9 \arity[u]_{#7} \arity "1";"3"^-{#2} } }
\newcommand{\functordn}[9]{
 \xymatrix{
    #4 \save[]+<0ex,5ex>*+{#1}="1"  \restore
       \arity@{}[rd]|{\longmapsto}
  & #5 \save[]+<0ex,5ex>*+{#3}="3"  \restore
      \arity[d]^{#7}
  \\
   #8  \arity[u]^{#6}  & #9 \arity "1";"3"^-{#2} } }

\newcommand{\larr}{->}
\newcommand{\rarr}{->}

\newcommand{\xfunctor}[9]{
 \xymatrix{
    #4 \save[]+<0ex,5ex>*+{#1}="1"  \restore
      \ifthenelse{\equal{\larr}{->}}{\arity[d]_{#6}}{}
      \ifthenelse{\equal{\larr}{<-}}{\arity[d];[]^{#6}}{}
      \ifthenelse{\equal{\larr}{-<}}{\arity@{< }[d]_{#6}}{}
      \arity@{}[rd]|{\longmapsto}
  & #5 \save[]+<0ex,5ex>*+{#3}="3"  \restore
      \ifthenelse{\equal{\rarr}{->}}{\arity[d]^{#7}}{}
      \ifthenelse{\equal{\rarr}{<-}}{\arity[d];[]_{#7}}{}
      \ifthenelse{\equal{\rarr}{-<}}{\arity@{< }[d]^{#7}}{}
  \\
   #8 & #9 \arity "1";"3"^-{#2} } }

\theoremstyle{plain}
\newtheorem{theorem}{Theorem}[section]
\newtheorem{proposition}[theorem]{Proposition}

\newtheorem{corollary}[theorem]{Corollary}

\newtheorem{remark}[theorem]{Remark}
\theoremstyle{definition}
\newtheorem{definition}[theorem]{Definition}

\newtheorem{examples}[theorem]{Examples}

\newtheorem*{assumption}{Assumption}

\newtheorem*{conventions}{Conventions}
\newtheorem*{convention}{Convention}

\theoremstyle{remark}

\newcommand{\arity}{\mathsf{ar}}

\newtheorem{claim}[theorem]{Claim}

\numberwithin{equation}{section}

\begin{document}
\title[Lallement functor]{Lallement functor is a weak right multiadjoint}

\author[Climent]{J. Climent Vidal}
\address{Universitat de Val\`{e}ncia\\
         Departament de L\`{o}gica i Filosofia de la Ci\`{e}ncia\\
         Av. Blasco Ib\'{a}\~{n}ez, 30-$7^{\mathrm{a}}$, 46010 Val\`{e}ncia, Spain}
\email{Juan.B.Climent@uv.es}
\author[Cosme]{E. Cosme Ll\'{o}pez}
\address{Universitat de Val\`{e}ncia\\
         Departament de Matem\`{a}tiques\\
         Dr. Moliner, 50, 46100 Burjassot, Val\`{e}ncia, Spain}
\email{enric.cosme@uv.es}

\keywords{ Grothendieck construction, Lallement systems of $\Sigma$-algebras, Lallement sum, naturally preordered idempotent $\Sigma$-algebras, functor of Lallement, semi-inductive Lallement systems of $\Sigma$-algebras, Lallement functor is a weak right multiadjoint, P{\l}onka functor.}
\date{May 24th, 2023}

\begin{abstract}
For a plural signature $\Sigma$ and with regard to the category $\mathsf{NPIAlg}(\Sigma)_{\mathsf{s}}$, of naturally preordered idempotent $\Sigma$-algebras and surjective homomorphisms, we define a contravariant functor $\mathrm{Lsys}_{\Sigma}$ from $\mathsf{NPIAlg}(\Sigma)_{\mathsf{s}}$ to 
$\mathsf{Cat}$, the category of categories, that assigns to $\mathbf{I}$ in 
$\mathsf{NPIAlg}(\Sigma)_{\mathsf{s}}$ the category $\mathbf{I}$-$\mathsf{LAlg}(\Sigma)$, of $\mathbf{I}$-semi-inductive Lallement systems of $\Sigma$-algebras, and a covariant functor $(\mathsf{Alg}(\Sigma)\,{\downarrow_{\mathsf{s}}}\, \cdot)$ from $\mathsf{NPIAlg}(\Sigma)_{\mathsf{s}}$ to $\mathsf{Cat}$, that assigns to $\mathbf{I}$ in $\mathsf{NPIAlg}(\Sigma)_{\mathsf{s}}$ the category $(\mathsf{Alg}(\Sigma)\,{\downarrow_{\mathsf{s}}}\, \mathbf{I})$, of the coverings of $\mathbf{I}$, i.e., the ordered pairs 
$(\mathbf{A},f)$ in which $\mathbf{A}$ is a $\Sigma$-algebra and $f\colon \mathbf{A}\mor \mathbf{I}$ a surjective homomorphism. Then, by means of the Grothendieck construction, we obtain the categories 
$\int^{\mathsf{NPIAlg}(\Sigma)_{\mathsf{s}}}\mathrm{Lsys}_{\Sigma}$ and
$\int_{\mathsf{NPIAlg}(\Sigma)_{\mathsf{s}}}(\mathsf{Alg}(\Sigma)\,{\downarrow_{\mathsf{s}}}\, \cdot)$;  define a functor $\mathfrak{L}_{\Sigma}$ from the first category to the second, which we will refer to as the Lallement functor; and prove that it is a weak right multiadjoint. Finally, we state the relationship between the P{\l}onka functor and the Lallement functor. 
\end{abstract}

\maketitle

\section{Introduction}

Let $\Sigma$ be a signature without $0$-ary operation symbols, $\mathbf{A}$ a $\Sigma$-algebra and 
$\Phi$ a congruence on $\mathbf{A}$. Then $\mathbf{A}/\Phi$ is an idempotent 
$\Sigma$-algebra if, and only if, for every $a\in A$, the equivalence class $[a]_{\Phi}$ of $a$ relative to $\Phi$ is a subalgebra of $\mathbf{A}$. Thus, a $\Sigma$-algebra $\mathbf{A}$ can be decomposed by congruences on $\mathbf{A}$ having the property that all equivalences classes are subalgebras of $\mathbf{A}$ if, and only if, it can be decomposed by congruences on $\mathbf{A}$ such that their quotients are idempotent $\Sigma$-algebras. Such a type of decomposition is referred to as an idempotent decomposition~\cite{CPB01}.

With respect to the just mentioned type of decompositions, and following the exposition in~\cite{CPB01}, two fundamental problems arise. The first of these has to do with how to break down an algebra into subalgebras in such a way that the structure of the algebra can be described in terms of the structure of its component subalgebras. Thus, e.g., if the congruence $\Phi$ on $\mathbf{A}$ is such that $\mathbf{A}/\Phi$ is an idempotent $\Sigma$-algebra, then we obtain a family 
$(\mathbf{A}_{[a]_{\Phi}})_{[a]_{\Phi}\in A/\Phi}$ of subalgebras of $\mathbf{A}$, indexed by the idempotent algebra $\mathbf{A}/\Phi$, where, for every $[a]_{\Phi}\in A/\Phi$, 
$\mathbf{A}_{[a]_{\Phi}}$ is the subalgebra of $\mathbf{A}$ whose underlying set is $\mathrm{pr}^{-1}_{\Phi}[\{[a]_{\Phi}\}] = [a]_{\Phi}$. Moreover, (1) 
$A = \bigcup_{[a]_{\Phi}\in A/\Phi} A_{[a]_{\Phi}}$, (2) for every $[a]_{\Phi}$, 
$[b]_{\Phi}\in A/\Phi$, if $[a]_{\Phi} \neq[b]_{\Phi}$, then 
$A_{[a]_{\Phi}}\cap A_{[a]_{\Phi}} = \varnothing$ and (3), for every $n\in \mathbb{N}-\{0\}$, every $\sigma\in \Sigma_{n}$ and every $(x_{i})_{i\in n}$, $(a_{i})_{i\in n}\in A^{n}$, if, for every $i\in n$, $(x_{i},a_{i})\in \Phi$, then $F^{\mathbf{A}}_{\sigma}((x_{i})_{i\in n})\in [F^{\mathbf{A}}_{\sigma}((a_{i})_{i\in n})]_{\Phi}$, i.e., 
$$
\textstyle
\{ F^{\mathbf{A}}_{\sigma}((x_{i})_{i\in n})\mid (x_{i})_{i\in n}\in \prod_{i\in n} A_{[a_{i}]_{\Phi}}\}\subseteq A_{[F^{\mathbf{A}}_{\sigma}((a_{i})_{i\in n})]_{\Phi}}.
$$

The second one, the composition problem, has to do with the following situation: given a family $(\mathbf{A}_{i})_{i\in I}$ of $\Sigma$-algebras indexed by the underlying set of an idempotent $\Sigma$-algebra $\mathbf{I}$ how to realize the operation symbols of $\Sigma$ as operations on $A=\coprod_{i\in I} A_{i}$ so that $A$ is equipped with a structure of $\Sigma$-algebra giving rise to a $\Sigma$-algebra $\mathbf{A}$; the equivalence relation on $A$ canonically associated to the partition $\{A_{i}\times \{i\}\mid i\in I\}$ of $A$ is a congruence relation on $\mathbf{A}$; and the associated quotient $\Sigma$-algebra is isomorphic to $\mathbf{I}$. 
 
From a category-theoretic standpoint, it could be said that the solution to the above problems, under certain conditions advanced in the abstract, is given, for the second, by the existence of the Lallement functor, whose object mapping assigns to a semi-inductive Lallement system of $\Sigma$-algebras its Lallement sum, and, for the first of them, by the existence of a weak left multiadjoint to the Lallement functor.

The source of these problems can be traced back to Clifford, who first investigated them in~\cite{C41}. Clifford's work focused on the study of semigroups that could be decomposed into semilattices of groups. He addressed the composition problem by employing inductive systems of homomorphisms over a semilattice, which we call  Clifford's construction, and here is the fundamental theorem, stated by Clifford in~\cite[Theorem 3, pp.~1044--1045]{C41}, in which such a construction is embodied:
\begin{quotation} 
Every semigroup which admits relative inverses and in which every pair of idempotent elements commute is isomorphic with a semigroup $S$ constructed as follows. 

Let $P$ be any semi-lattice, and to each $\alpha$ in $P$ assign a group $S_{\alpha}$, such that no two of them have an element in common. To each pair of elements $\alpha>\beta$ of $P$ assign a homomorphism $\phi_{\alpha\beta}$ of $S_{\alpha}$ into $S_{\beta}$ such that if $\alpha > \beta > \gamma$ then 
$$
\phi_{\alpha\beta}\phi_{\beta\gamma} = \phi_{\alpha\gamma}.
$$ 
Let $\phi_{\alpha\alpha}$ be the identical automorphism of $S_{\alpha}$. Let $S$ be the class sum of the groups $S_{\alpha}$, and define the product of any two elements $a_{\alpha}$, $b_{\beta}$ of $S$ ($a_{\alpha}$ in $S_{\alpha}$ and $b_{\beta}$ in $S_{\beta}$) by 
$$
a_{\alpha}b_{\beta} = (a_{\alpha}\phi_{\alpha\gamma})(b_{\beta}\phi_{\beta\gamma})
$$
where $\gamma = \alpha\beta$ is the product of $\alpha$ and $\beta$ in $P$.

Conversely, any semigroup $S$ constructed in this fashion admits relative inverses, and every idempotent element of $S$ is in the center of $S$. 

The semi-lattice $P$, the groups $S_{\alpha}$ and the homomorphisms $\phi_{\alpha\beta}$ together constitute a complete set of invariants of $S$.
\end{quotation}
The work of Clifford was generalized in several directions. Thus, 
P{\l}onka, in~\cite{P67}, introduced a new construction in universal algebra, which is  closely related to the classical notion of inductive limit of an inductive system of algebras relative to (or indexed by) a directed preorder (see~\cite{CC24} for a category-theoretic explanation of this connection). Specifically, for a sup-semilattice inductive system of algebras, he defined its sum, later called the P{\l}onka sum, and proved, by means of the notion of partition function of the underlying set of an algebra without $0$-ary operations, that the correspondence between partition functions and representations as the sum of a sup-semilattice inductive system of algebras is one-to-one. 
Another generalization of Clifford's construction (presented, as mentioned above,  in~\cite[Theorem 3, pp.~1044--1045]{C41} and repeated in~\cite[Theorem 4.11, p.~128]{CP61}) was given by Lallement in~\cite[Th\'{e}or\`{e}me 2.19, pp.~83--84]{La67}. He defined the sum of a semilattice ordered system of semigroups by using a not necessarily transitive system of homomorphisms over a semilattice. Here is the statement of the theorem we have just mentioned: 
\begin{quotation}
Let $Y$ be a semilattice. To each $\alpha\in Y$ is assigned a completely simple semigroup $S_{\alpha}$ so that 
$S_{\alpha}\cap S_{\beta} = \varnothing$ if $\alpha \neq \beta$. To each pair 
$(\alpha,\beta)\in Y\times Y$ such that $\alpha\geq\beta$ is assigned a mapping 
$\Phi_{\alpha,\beta}$ from $S_{\alpha}$ to the translational hull 
$\overline{S_{\beta}}$ of $S_{\beta}$, the mappings $\Phi_{\alpha,\beta}$ being subject to the following conditions
\begin{enumerate}
\item[1)] $S_{\alpha}\Phi_{\alpha,\alpha\beta}\cdot S_{\beta}\Phi_{\beta,\alpha\beta}\subseteq S_{\alpha\beta}\Phi_{\alpha\beta,\alpha\beta}$, for every 
$\alpha,\beta\in Y$;
\item[2)] for every $\alpha$, $\beta$, $\gamma\in Y$ such that $\alpha\beta\geq \gamma$ and for every $a_{\alpha}\in S_{\alpha}$ and $b_{\beta}\in S_{\beta}$: 
$$
[a_{\alpha}\Phi_{\alpha,\alpha\beta}\cdot b_{\beta}\Phi_{\beta,\alpha\beta}]\Phi_{\alpha\beta,\gamma} = a_{\alpha}\Phi_{\alpha,\gamma}\cdot b_{\beta}\Phi_{\beta,\gamma};
$$
\item[3)] $\Phi_{\alpha,\alpha}$ is the canonical inclusion of $S_{\alpha}$ into 
$\overline{S_{\alpha}}$.
\end{enumerate} 

On $D = \bigcup_{\alpha\in Y}S_{\alpha}$ an operation is defined by $a_{\alpha}\bigcdot b_{\beta} = a_{\alpha}\Phi_{\alpha,\alpha\beta}\cdot b_{\beta}\Phi_{\beta,\alpha\beta}$ [Op.]. Then $D$ is a semigoup union of groups. Reciprocally, every semigroup union of groups can be obtained in this way. 
\end{quotation}
Note that if in the third condition we identify $S_{\alpha}$ with its image, then the first condition can, simply, be written as $\Phi_{\alpha,\alpha\beta}[S_{\alpha}]\cdot\Phi_{\beta,\alpha\beta}[S_{\beta}]\subseteq S_{\alpha\beta}$, where $\Phi_{\alpha,\alpha\beta}[S_{\alpha}]$ is the direct image of $S_{\alpha}$ under $\Phi_{\alpha,\alpha\beta}$ and the same for $\Phi_{\beta,\alpha\beta}[S_{\beta}]$. 
Later on, Romanowska and Smith in~\cite{RS91}, defined, for Lallement systems of algebras, the corresponding sums, which they called Lallement sums. They also proved that every algebra decomposed into a semilattice sum of algebras can be constructed using this method. The reader can find more information about Lallement and  P{\l}onka sums in the survey paper~\cite{PR92} by P{\l}onka and Romanowska.

Lallement sums of direct systems of algebras are applied to several algebraic theories, including semigroup theory and semiring theory. These investigations were further pursued, among others, by Romanowska and Smith in~\cite{RS85, RS91, RS02}.

The main goal of this paper is to conduct a category-theoretic analysis of the Lallement sum, which we also refer to as the Lallement construction. Specifically, we show (1) that the Lallement construction is a functor, we will refer to it as the functor of Lallement, between categories obtained by the Grothendieck construction applied to a contravariant functor, for the domain of the functor, and to a covariant functor, for the codomain of the functor; and (2) that the functor of Lallement is a weak right multiadjoint. The just mentioned property is a generalization of the one of right multiadjoint originally introduced by Diers, among other places, in~\cite{YD84}.

The paper is structured as follows. In Section 2, in order to keep the exposition as self-contained as possible, we present those fundamental results in universal algebra, lattice theory, semigroup theory and category theory that will be used in the paper. We, in particular, recall the concepts of split fibration, split opfibration, split indexed category and opposite split indexed category, along with their associated Grothendieck constructions. We conclude this section by defining the concept of weak right multiadjoint functor. This will be the key concept of the paper, since we will show that the Lallement functor, once defined, falls under it.

In Section 3, for a $\Sigma$-algebra $\mathbf{I}$, we introduce the concept of a Lallement system of $\Sigma$-algebras relative to $\mathbf{I}$ and the notion of a morphism between Lallement systems of $\Sigma$-algebras relative to $\mathbf{I}$. This leads us to define an associated category. Moreover, we define a contravariant functor $\mathrm{Lsys}_{\Sigma}$ from the category $\mathsf{Alg}(\Sigma)$ to the category $\mathsf{Cat}$ and from it we obtain, by means of the Grothendieck construction, the category $\int^{\mathsf{Alg}(\Sigma)}\mathrm{Lsys}_{\Sigma}$, of labeled Lallement systems of $\Sigma$-algebras.

In Section 4, after defining the notion of plural signature, we define the notion of idempotent $\Sigma$-algebra and provide a characterization of the idempotency of a $\Sigma$-algebra in terms of constant mappings. After that, we define the concept of naturally preordered $\Sigma$-algebra,  
the category of naturally preordered idempotent  $\Sigma$-algebras and provide several characterizations of the concept of naturally preordered idempotent  $\Sigma$-algebra. 
Moreover, from the category $\mathsf{NPIAlg}(\Sigma)_{\mathsf{s}}$, of naturally preordered idempotent $\Sigma$-algebras and surjective homomorphisms, we obtain, on the one hand, by restricting
the contravariant functor $\mathrm{Lsys}_{\Sigma}$ to it the category $\int^{\mathsf{NPIAlg}(\Sigma)_{\mathsf{s}}}\mathrm{Lsys}_{\Sigma}$ which has as objects the Lallement systems of $\Sigma$-algebras labeled by naturally preordered idempotent $\Sigma$-algebras; and, on the other hand, after defining another functor 
$(\mathsf{Alg}(\Sigma)\,{\downarrow_{\mathsf{s}}}\, \bigcdot)$ from $\mathsf{NPIAlg}(\Sigma)_{\mathsf{s}}$ to $\mathsf{Cat}$, the category $\int_{\mathsf{NPIAlg}(\Sigma)_{\mathsf{s}}}(\mathsf{Alg}(\Sigma)\,{\downarrow_{\mathsf{s}}}\, \bigcdot)$ which has as objects the ordered pairs 
$(\mathbf{I},(\mathbf{A},f))$, where $\mathbf{I}$ is a naturally preordered idempotent $\Sigma$-algebra and $(\mathbf{A},f)$ an object of $(\mathsf{Alg}(\Sigma)\,{\downarrow_{\mathsf{s}}}\, \mathbf{I})$, i.e., a $\Sigma$-algebra 
$\mathbf{A}$ and a surjective homomorphism $f$ from $\mathbf{A}$ to $\mathbf{I}$. 

In Section 5, for a plural signature $\Sigma$, we introduce the concept of the Lallement sum for a Lallement system of $\Sigma$-algebras relative to a naturally preordered idempotent $\Sigma$-algebra. This enables us to define a covariant functor, denoted as $\mathfrak{L}_{\Sigma}$, from $\int^{\mathsf{NPIAlg}(\Sigma)_{\mathsf{s}}}\mathrm{Lsys}_{\Sigma}$ to $\int_{\mathsf{NPIAlg}(\Sigma)_{\mathsf{s}}}(\mathsf{Alg}(\Sigma),{\downarrow_{\mathsf{s}}}, \bigcdot)$, which we shall refer to as the Lallement functor relative to $\Sigma$.

In Section 6, after introducing the notion of semi-inductive Lallement systems of $\Sigma$-algebras relative to a naturally preordered idempotent $\Sigma$-algebra, we prove that the restriction of the functor $\mathfrak{L}_{\Sigma}$ to these Lallement systems is a weak right multiadjoint, and this is the category-theoretical content of the Lallement construction.

Finally, in Section 7, after proving that there exists a functor $\int^{G}\eta$ from $\int^{\mathsf{Sl}_{\mathsf{s}}}\mathrm{Isys_{\Sigma}}$ to $\int^{\mathsf{NPIAlg}(\Sigma)_{\mathsf{s}}}\mathrm{Lsys}_{\Sigma}$; defining the category 
$\mbox{\sffamily{\upshape{P{\l}Alg}}}(\Sigma)$, of P{\l}onka algebras, and its  subcategory $\mbox{\sffamily{\upshape{P{\l}Alg}}}(\Sigma)_{\mathsf{s}}$, consisting of all objects of $\mbox{\sffamily{\upshape{P{\l}Alg}}}(\Sigma)$ and with morphisms the surjective homomorphisms between P{\l}onka algebras; defining functors $P$ from 
$\mbox{\sffamily{\upshape{P{\l}Alg}}}(\Sigma)_{\mathsf{s}}$ to $\mathsf{Alg}(\Sigma)_{\mathsf{s}}$, the category of $\Sigma$-algebras and surjective homomorphisms, and $Q$ from $\int_{\mathsf{NPIAlg}(\Sigma)_{\mathsf{s}}}(\mathsf{Alg}(\Sigma)\,{\downarrow_{\mathsf{s}}}\, \bigcdot)$ to  
$\mathsf{Alg}(\Sigma)_{\mathsf{s}}$; and recalling the essentials about $\mbox{\upshape{P{\l}}}_{\Sigma}$, the P{\l}onka functor from 
$\int^{\mathsf{Sl}}\mathrm{Isys}_{\Sigma}$ to $\mbox{\sffamily{\upshape{P{\l}Alg}}}(\Sigma)$ (for full details about this functor see~\cite{CC24}), we state the relationship between the birestriction of $\mbox{\upshape{P{\l}}}_{\Sigma}$ to $\int^{\mathsf{Sl}_{\mathsf{s}}}\mathrm{Isys}_{\Sigma}$ and $\mbox{\sffamily{\upshape{P{\l}Alg}}}(\Sigma)_{\mathsf{s}}$, also denoted by $\mbox{\upshape{P{\l}}}_{\Sigma}$, and the Lallement functor. Concretely, we prove that there exists a natural isomorphism $\theta$ from $P\circ \mbox{\upshape{P{\l}}}_{\Sigma}$ to $Q\circ\mathfrak{L}_{\Sigma}\circ\int^{G}\eta$. 

\begin{conventions}
Our underlying set theory is $\mathbf{ZFSk}$, Zermelo-Fraenkel-Skolem set theory (also known as $\mathbf{ZFC}$, i.e., Zermelo-Fraenkel set theory with the axiom of choice) plus the existence of a Grothendieck universe $\mathbf{U}$, fixed once and for all (see~\cite[pp.~21--24]{sM98}). The $\mathbf{U}$-small sets are those sets which are elements of $\mathbf{U}$ and the $\mathbf{U}$-large sets are those sets which are subsets of $\mathbf{U}$.

In all that follows we use standard concepts and constructions from set theory, see, e.g.,~\cite{nB70, end77}; universal algebra, see, e.g.,~\cite{gb15, bl05, bs81, gG08, w92}; lattice theory, see, e.g.,~\cite{bl05, HK71}; semigroup theory, see, e.g.,~\cite{JH80}; and category theory, see, e.g.,~\cite{BW85, BW12, Gro71, hs73, sM98, em76}. Nevertheless, regarding set theory, we have adopted the following conventions. An \emph{ordinal} $\alpha$ is a transitive set that is well-ordered by $\in$; thus, $\alpha=\{\beta \mid \beta \in \alpha\}$. The first transfinite ordinal $\omega_{0}$ will be denoted by $\mathbb{N}$, which is the set of all \emph{natural numbers}, and, from what we have just said about the ordinals, for every $n\in\mathbb{N}$, $n=\{0,\dots, n-1\}$. For a set $A$ we denote by 
$\mathrm{Sub}(A)$ the set of all subsets of $A$. A function from $A$ to $B$ is a subset $F$ of $A\times B$ satisfying the functional condition and a mapping from $A$ to $B$ is an ordered triple $f=(A,F,B)$, denoted by $f\colon A\mor B$, in which $F$ is a function from $A$ to $B$. 
If $f$ is a mapping from $A$ to $B$, then the mapping $f[\cdot]$ from $\mathrm{Sub}(A)$ to $\mathrm{Sub}(B)$, of $f$-direct image formation, sends $X$ in $\mathrm{Sub}(A)$ to
$
f[X] = \{y\in B\mid \exists\,x\in X\,(y = f(x))\}
$
in $\mathrm{Sub}(B)$. Finally, we let $\mathrm{Ker}(f)$ stand for the kernel of $f$ and $\mathrm{Im}(f)$ for the image of $f$, i.e., for $f[A]$.
Moreover, with regard to category theory, for a category $\mathsf{C}$ we denote by 
$\mathrm{Ob}(\mathsf{C})$ the set of objects of $\mathsf{C}$, by $\mathrm{Mor}(\mathsf{C})$ the set of morphisms of $\mathsf{C}$ and, for $x$, $y\in \mathrm{Ob}(\mathsf{C})$, we denote by $\mathrm{Hom}(x,y)$ the set of all morphisms in $\mathsf{C}$ from $x$ to $y$. Moreover, we let $\mathrm{sc}_{\mathsf{C}}$ and $\mathrm{tg}_{\mathsf{C}}$ stand for the mappings from $\mathrm{Mor}(\mathsf{C})$ to $\mathrm{Ob}(\mathsf{C})$ that assign to a morphism of $\mathsf{C}$ its source and its target, respectively; and for a functor $F$ from a category $\mathsf{C}$ to another $\mathsf{D}$ we sometimes speak of the object mapping of $F$, i.e., of the mapping from $\mathrm{Ob}(\mathsf{C})$ to $\mathrm{Ob}(\mathsf{D})$ that sends an object $x$ of $\mathsf{C}$ to the object $F(x)$ of $\mathsf{D}$. In addition, for a category 
$\mathsf{C}$ we let $\mathsf{C}^{\mathrm{op}}$ stand for the dual (or opposite) of $\mathsf{C}$ and $F^{\mathrm{op}}\colon \mathsf{C}^{\mathrm{op}}\mor \mathsf{D}^{\mathrm{op}}$ for the dual (or opposite) of $F\colon \mathsf{C}\mor \mathsf{D}$. Finally, we let $\mathsf{Set}$ stand for the category of sets in $\mathbf{U}$ and mappings between them and $\mathsf{Cat}$ for the category of all $\mathbf{U}$-categories, as defined in~\cite{GW23}, p.~738, i.e., the category whose objects are all those categories $\mathsf{C}$ such that $\mathrm{Ob}(\mathsf{C})$ is a subset of $\mathbf{U}$ and, for every $x$, $y\in \mathrm{Ob}(\mathsf{C})$, $\mathrm{Hom}(x,y)\in \mathbf{U}$ (in other words, the category whose objects are both  $\mathbf{U}$-large categories, i.e., those for which $\mathrm{Ob}(\mathsf{C})$ and $\mathrm{Mor}(\mathsf{C})$ are $\mathbf{U}$-large sets, and $\mathbf{U}$-locally small, i.e., those for which $\mathrm{Hom}(x,y)$ is a $\mathbf{U}$-small set, for every $x$, $y\in \mathrm{Ob}(\mathsf{C})$); and whose morphisms are the functors between $\mathbf{U}$-categories. 
\end{conventions}

\section{Preliminaries}

In this section we collect those basic facts about universal algebra, lattice theory, semigroup theory  and category theory, that we will need to carry out our research. For a full treatment of the topics of this section we refer the reader to \cite{BW85, gb15, bl05, bs81, gG08, Gro71, hs73, em76, w92}.

We begin by defining the concept of signature.

\begin{definition}\label{DSig}
A \emph{signature} is a mapping $\Sigma = (\Sigma_{n})_{n\in\mathbb{N}}$ from $\mathbb{N}$ to 
$\mathbf{U}$ which sends a natural number $n\in\mathbb{N}$ to the set $\Sigma_{n}$ of the \emph{$n$-ary formal operations} or \emph{$n$-ary operation symbols}. 
\end{definition}

\begin{assumption}
From now on, unless otherwise stated, $\Sigma$ stands for a signature as set forth in Definition~\ref{DSig}, fixed once and for all.
\end{assumption}

We shall now define the concepts of algebra and homomorphism between algebras.

\begin{definition}\label{DAlg}
A $\Sigma$-\emph{algebra} is a pair $(A,F)$, denoted by $\mathbf{A}$, where $A$ is a set and $F$ a structure of $\Sigma$-algebra on $A$, i.e., a family $(F_{n})_{n\in\mathbb{N}}$ where, for $n\in\mathbb{N}$, $F_{n}$ is a mapping from $\Sigma_{n}$ to $\mathrm{Hom}(A^{n},A)$ such that, for every $(\sigma,n)\in\coprod_{n\in\mathbb{N}} \Sigma_{n}$, $F_{n}(\sigma)$, written as $F_{\sigma}$ for short, is a mapping from $A^{n}$ to $A$. In some cases, to avoid any confusion, we will denote by $F^{\mathbf{A}}$ the structure of $\Sigma$-algebra on $A$, and, for $n\in \mathbb{N}$ and $\sigma\in \Sigma_{n}$, by $F^{\mathbf{A}}_{\sigma}$ the corresponding operation. 
A \emph{homomorphism} from $\mathbf{A}$ to $\mathbf{B}$, where $\mathbf{B} = (B,G)$, is a triple $(\mathbf{A},f,\mathbf{B})$, denoted by $f\colon \mathbf{A}\mor \mathbf{B}$, where $f$ is a  mapping from $A$ to $B$ such that, for every $n\in\mathbb{N}$, every  $\sigma\in \Sigma_{n}$ and every $(a_{j})_{j\in n}\in A^{n}$, we have 
$$
f(F^{\mathbf{A}}_{\sigma}((a_{j})_{j\in n})) = G^{\mathbf{B}}_{\sigma}((f(a_{j}))_{j\in n}).
$$
We will denote by $\mathsf{Alg}(\Sigma)$ the category of $\Sigma$-algebras and homomorphisms and by $\mathrm{Alg}(\Sigma)$ the set of objects  of $\mathsf{Alg}(\Sigma)$.

\end{definition}

\begin{assumption}
From now on, we will assume that the underlying set of any $\Sigma$-algebra, semigroup or semilattice occurring in this work is an element of $\mathbf{U}$.
\end{assumption}

We shall now give the definitions of the concepts of subalgebra and of closed subset  of an algebra.

\begin{definition}\label{DSubAlg}
Let $\mathbf{A} = (A,F)$ be a $\Sigma$-algebra. We will say that a $\Sigma$-algebra 
$\mathbf{X}$ is a \emph{subalgebra} of $\mathbf{A}$ if $X\subseteq A$ and $\mathrm{in}_{X,A}$, the canonical embedding of $X$ into $A$, determines an injective homomorphism 
$\mathrm{in}_{\mathbf{X},\mathbf{A}} = (\mathbf{X},\mathrm{in}_{X,A},\mathbf{A})$ from $\mathbf{X}$ to $\mathbf{A}$. Let $X$ be a subset of $A$. We will say that $X$ is a \emph{closed subset} of $\mathbf{A}$ if $X$ is closed under the operations of 
$\mathbf{A}$, i.e., if, for every $n\in \mathbb{N}$, every $\sigma\in\Sigma_{n}$ and every $(a_{j})_{j\in n}\in X^{n}$, $F^{\mathbf{A}}_{\sigma}((a_{j})_{j\in n})\in X$. If $X$ is a closed subset of $\mathbf{A}$, then we denote by $\mathbf{X}$ the subalgebra of $\mathbf{A}$ canonically associated to $X$.

\end{definition}

We next give the definitions of the concepts of congruence on an algebra and of quotient of an algebra by a congruence on it.

\begin{definition}\label{DCong}
Let $\mathbf{A}$ be a $\Sigma$-algebra and $\Phi$ an equivalence on $A$. We will say that $\Phi$ is a \emph{congruence on} $\mathbf{A}$ if, for every $n\in \mathbb{N}-\{0\}$, every $\sigma\in \Sigma_{n}$,
and every $(a_{j})_{j\in n}, (b_{j})_{j\in n}\in A^{n}$, if, for every $j\in  n$, $(a_{j}, b_{j})\in\Phi$, then $$(F^{\mathbf{A}}_{\sigma}((a_{j})_{j\in n}), F^{\mathbf{A}}_{\sigma}((b_{j})_{j\in n}))\in \Phi.$$
We will denote by $\mathrm{Cgr}(\mathbf{A})$ the set of all congruences on $\mathbf{A}$.

%
For a congruence $\Phi$ on $\mathbf{A}$, the \emph{quotient $\Sigma$-algebra of} $\mathbf{A}$ \emph{by} $\Phi$, denoted by $\mathbf{A}/\Phi$, is the $\Sigma$-algebra $(A/\Phi, F^{\mathbf{A}/\Phi})$, where, for every  $n\in \mathbb{N}$ and every $\sigma\in \Sigma_{n}$, the operation $F_{\sigma}^{\mathbf{A}/\Phi}$ from $(A/\Phi)^{n}$ to $A/\Phi$, sends $([a_{j}]_{\Phi})_{j\in n}$ in $(A/\Phi)^{n}$ to $[F^{\mathbf{A}}_{\sigma}((a_{j})_{j\in n})]_{\Phi}$ in $A/\Phi$,
and the \emph{canonical projection} from $\mathbf{A}$ to $\mathbf{A}/\Phi$, denoted by $\mathrm{pr}_{\Phi}\colon \mathbf{A}\mor \mathbf{A}/\Phi$, is the surjective homomorphism determined by the projection from $A$ to $A/\Phi$.
\end{definition}

We next define the algebraic predecessor relation on an algebra and the reflexive and transitive closure of it. We recall that these notions were defined first by L\"{o}wig in~\cite{HL52} and were thoroughy investigated by Diener in~\cite{KD66, KD93}, who, moreover, in~\cite{KD93}, provides a careful and detailed historical study of them.

\begin{definition}\label{DSubPOrd} 
Let $\mathbf{A}=(A,F)$ be a $\Sigma$-algebra. We let $\prec_{\mathbf{A}}$ denote the binary relation on $A$ consisting of the ordered pairs $(b,a)\in A\times A$ for which there exists a natural number $n\in\mathbb{N}-\{0\}$, an operation symbol $\sigma\in \Sigma_{n}$ and a family of elements $(a_{j})_{j\in n}\in A^{n}$ such that $a=F^{\mathbf{A}}_{\sigma}((a_{j})_{j\in n})$ and, for some $j\in n$, $b=a_{j}$. If $(b,a)\in\prec_{\mathbf{A}}$, also written as $b\prec_{\mathbf{A}} a$, then we will say that $b$ \emph{$\prec_{\mathbf{A}}$-precedes} $a$ or, when no confusion can arise, that $b$ precedes $a$. We will call $\prec_{\mathbf{A}}$ the \emph{algebraic predecessor relation} on $\mathbf{A}$. We will denote by $\leq_{\mathbf{A}}$ the reflexive and transitive closure of $\prec_{\mathbf{A}}$, i.e., the preorder on $A$ generated by $\prec_{\mathbf{A}}$, and by $<_{\mathbf{A}}$ the transitive closure of $\prec_{\mathbf{A}}$. If $(b,a)\in \leq_{\mathbf{A}}$, also written as $b\leq_{\mathbf{A}} a$, then we will say that $b$ is an \emph{algebraic predecessor of $a$}. The preorder $\leq_{\mathbf{A}}$ will be called the \emph{algebraic predecessor preorder} on $\mathbf{A}$.
\end{definition}


\begin{remark}\label{RSubPOrd}
The preorder $\leq_{\mathbf{A}}$ on $A$ is $\bigcup_{n\in\mathbb{N}}\prec_{\mathbf{A}}^{n}$, where $\prec_{\mathbf{A}}^{0}$ is the diagonal of $A$ and, for $n\in\mathbb{N}$, $\prec_{\mathbf{A}}^{n+1}=\prec^{n}_{\mathbf{A}}\circ \prec_{\mathbf{A}}$. Thus, for every $(b,a)\in A\times A$, $b\leq_{\mathbf{A}} a$ if and only if $b=a$ or there exists a natural number $m\in\mathbb{N}-\{0\}$, a family of elements $(r_{k})_{k\in m+1}$ in $A^{m+1}$ such that $r_{0}=b$, $r_{m}=a$ and, for every $k\in m$, $r_{k}\prec_{\mathbf{A}} r_{k+1}$. Moreover, $<_{\mathbf{A}}$, the transitive closure of $\prec_{\mathbf{A}}$ is $\bigcup_{n\in\mathbb{N}-\{0\}}\prec_{\mathbf{A}}^{n}$.
\end{remark}

\begin{proposition}\label{PSubOrd}
Let $f$ be a homomorphism from a $\Sigma$-algebra $\mathbf{A}$ to another $\Sigma$-algebra $\mathbf{B}$. Then, for every $x$, $y\in A$, if $y\leq_{\mathbf{A}}x$, then 
$f(y)\leq_{\mathbf{B}}f(x)$.
\end{proposition}


We shall now recall the notions of semilattice and of idempotent commutative semigroup and state the well known fact that the corresponding categories are isomorphic.

\begin{definition}\label{DSup}
A \emph{semilattice} is an ordered set $\mathbf{I} = (I,\leq)$ in which any two elements $i$, $j$ have a least upper bound $i\vee j$ (equivalent terminology for this is a \emph{sup-semilattice}, \emph{join semilattice} or a $\vee$-\emph{semilattice}). A morphism from a semilattice $\mathbf{I}$ to another semilattice $\mathbf{P}$ is a triple $(\mathbf{I},\xi,\mathbf{P})$, abbreviated to $\xi\colon \mathbf{I}\mor \mathbf{P}$, where $\xi$ is a mapping from $I$ to $P$ such that, for every $i$, $j\in I$, $\xi(i\vee j) = \xi(i)\vee \xi(j)$. We let $\mathsf{Sl}$ stand for the category of \emph{semilattices} and morphisms between them.
\end{definition}

\begin{remark}
If $\xi\colon \mathbf{I}\mor \mathbf{P}$ is a morphism between semilattices, then $\xi$ is an isotone mapping between them. Moreover, in $\mathsf{Sl}$, the monomorphisms are the injective morphisms and the epimorphisms the surjective morphisms.
\end{remark}


\begin{definition}\label{ICSgr}
A \emph{semigroup} is an ordered pair $\mathbf{A} = (A,\cdot)$ in which $\cdot$ is a binary operation on $A$ such that, for every $x$, $y$, $z\in A$, $(x\cdot y)\cdot z = x\cdot (y\cdot z)$. A homomorphism from a semigroup $\mathbf{A}$ to another semigroup $\mathbf{B}$ is a triple 
$(\mathbf{A},f,\mathbf{B})$, abbreviated to $f\colon \mathbf{A}\mor \mathbf{B}$, where $f$ is a mapping from $A$ to $B$ such that, for every $x$, $y\in A$, $f(x\cdot y) = f(x)\cdot f(y)$. A semigroup $\mathbf{A}$ is \emph{commutative} if, for every $x$, $y\in A$, $x\cdot y = y\cdot x$.
A semigroup $\mathbf{A}$ is \emph{idempotent} if, for every $x\in A$, $x\cdot x = x$. We let 
$\mathsf{Sgr}$, $\mathsf{CSgr}$ and $\mathsf{ICSgr}$ stand for the categories of \emph{semigroups}, \emph{commutative semigroups} and \emph{idempotent commutative semigroups}, respectively. 
\end{definition}

\begin{proposition}
The category $\mathsf{Sl}$ is isomorphic to the category $\mathsf{ICSgr}$.
\end{proposition}

In what follows we begin by defining the notions of cartesian morphism, fibration, morphism between fibrations, split fibration, morphism between split fibrations as well as the notions of opfibration and split opfibration---and in doing so we follow the presentation by Barr and Wells in~\cite{BW12}, pp.~327--329, and by Jacobs in~\cite{bj99}, p.~73. After this we define the notions of split $\mathsf{I}$-indexed category, for a category $\mathsf{I}$, split indexed category, morphism between split indexed categories as well as the notions of opposite split $\mathsf{I}$-indexed category, and opposite split indexed category---also following the terminology used by Jacobs in~\cite{bj99}, p.~51. Then we define the Grothendieck construction (see~\cite{Gro71}) by which we mean both the functor from the category of split indexed categories to the category of split fibrations and the functor from the category of opposite split indexed categories to the category of split opfibrations.

\begin{definition}
Let $\pi\colon \mathsf{E}\mor \mathsf{B}$ be a functor, $f\colon c\mor d$ a morphism of $\mathsf{B}$ and $y$ an object of $\mathsf{E}$ such that $\pi(y) = d$. We will say that a morphism $u\colon x\mor y$ of $\mathsf{E}$ is \emph{cartesian} for $f$ and $y$ if $\pi(u) = f$ and, for every object $z$ of $\mathsf{E}$, every morphism $v\colon z\mor y$ of $\mathsf{E}$ and every morphism $h\colon \pi(z)\mor c$ of $\mathsf{B}$ for which $f\circ h = \pi(v)$, there exists a unique morphism $w\colon z\mor x$ of $\mathsf{E}$ such that $u\circ w = v$ and $\pi(w) = h$. We let $\mathrm{Cart}(f,y)$ stand for the set of all cartesian morphisms for $f$ and $y$.

A functor $\pi\colon \mathsf{E}\mor \mathsf{B}$ is a \emph{fibration} if there exists a cartesian morphism for every morphism $f\colon c\mor d$ of $\mathsf{B}$ and every object $y$ of $\mathsf{E}$ for which $\pi(y) = d$. Sometimes a fibration $\pi\colon \mathsf{E}\mor \mathsf{B}$ is called a \emph{fibred category} or a \emph{category fibred over} $\mathsf{B}$ and denoted by $(\mathsf{E},\pi)$. A \emph{morphism of fibrations} from a fibration $\pi\colon \mathsf{E}\mor \mathsf{B}$ to another $\pi'\colon \mathsf{E}'\mor \mathsf{B}'$ is a pair of functors $H\colon \mathsf{E}\mor \mathsf{E}'$ and $K\colon \mathsf{B}\mor \mathsf{B}'$ such that $K\circ \pi = \pi'\circ H$ and $H$ sends cartesian morphisms in $\mathsf{E}$ to cartesian morphisms in $\mathsf{E}'$. We let $\mathsf{Fib}$ stand for the corresponding category.

A functor $\pi\colon \mathsf{E}\mor \mathsf{B}$ is an \emph{opfibration} if 
$\pi^{\mathrm{op}}\colon \mathsf{E}^{\mathrm{op}}\mor \mathsf{B}^{\mathrm{op}}$ is a fibration. We let $\mathsf{OpFib}$ stand for the corresponding category. For details see~\cite{BW12}.

Let $\mathrm{Mor}(\mathsf{B})\times_{\mathrm{Ob}(\mathsf{B})}\mathrm{Ob}(\mathsf{E})$ be the fibered product of $\mathrm{tg}_{\mathsf{B}}\colon \mathrm{Mor}(\mathsf{B})\mor \mathrm{Ob}(\mathsf{B})$, the mapping that sends a morphism of $\mathsf{B}$ to its target, and $\pi\colon \mathrm{Ob}(\mathsf{E})\mor \mathrm{Ob}(\mathsf{B})$, the object mapping of the functor $\pi$. Then a \emph{clivage} for a fibration $\pi\colon \mathsf{E}\mor \mathsf{B}$ is an element $\gamma$ of $\prod_{(f,y)\in \mathrm{Mor}(\mathsf{B})\times_{\mathrm{Ob}(\mathsf{B})}\mathrm{Ob}(\mathsf{E})}\mathrm{Cart}(f,y)$, i.e., a mapping that assigns to every pair $(f,y)$, in which $f\colon c\mor d$ is a morphism of $\mathsf{B}$ and $y$ an object of $\mathsf{E}$ for which $\pi(y) = d$, a morphism $\gamma_{f,y}$ of $\mathsf{E}$ that is cartesian for $f$ and $y$. The clivage  $\gamma$ is a \emph{splitting} of the fibration $\pi\colon \mathsf{E}\mor \mathsf{B}$ if it is such that (1) for every object $d$ of $\mathsf{B}$ and every object $y$ of $\mathsf{E}$ for which $\pi(y) = d$, then $\gamma_{\mathrm{id}_{d},y} = \mathrm{id}_{y}$ and (2)  for every pair of morphisms $f\colon c\mor d$ and $g\colon d\mor e$ of 
$\mathsf{B}$ and every pair of objects $y$ and $z$ of $\mathsf{E}$ if $\pi(z) = e$ and $y$ is the domain of $\gamma_{g,z}$, then $\gamma_{g,z}\circ\gamma_{f,y} = \gamma_{g\circ f,z}$.

A fibration $\pi\colon \mathsf{E}\mor \mathsf{B}$ is \emph{split} if it has a splitting. A \emph{morphism of split fibrations} from a split fibration $\pi\colon \mathsf{E}\mor \mathsf{B}$ to another $\pi'\colon \mathsf{E}'\mor \mathsf{B}'$ is a pair of functors $(H,K)$ as above where $H$ preserves the splitting up to equality (not up to isomorphism). We let $\mathsf{SFib}$ stand for the corresponding category.

An opfibration $\pi\colon \mathsf{E}\mor \mathsf{B}$ is \emph{split} if 
$\pi^{\mathrm{op}}\colon \mathsf{E}^{\mathrm{op}}\mor \mathsf{B}^{\mathrm{op}}$ has a splitting.  We let $\mathsf{SOpFib}$ stand for the corresponding category. For details see~\cite{BW12}.
\end{definition}

\begin{definition}
Let $\mathsf{I}$ be a category. A \emph{split} $\mathsf{I}$-\emph{indexed category} is a contravariant functor $F$ from $\mathsf{I}$ to $\mathsf{Cat}$. Given an object $i\in \mathrm{Ob}(\mathsf{I})$, we write $\mathsf{F}_{i}$ for the category $F(i)$, and given a morphism $\varphi\in \mathrm{Hom}_{\mathsf{I}}(i,j)$, we write $F_{\varphi}$ for the functor $F(\varphi)\colon \mathsf{F}_{j}\mor \mathsf{F}_{i}$. A \emph{split indexed category} is an ordered pair $(\mathsf{I},F)$ in which $\mathsf{I}$ is a category and $F$ a split $\mathsf{I}$-indexed category. A \emph{morphism of split indexed categories} from a split indexed category $(\mathsf{I},F)$ to another $(\mathsf{I}',F')$ is an ordered pair $(G,\eta)$, where $G$ is a functor from $\mathsf{I}$ to $\mathsf{I}'$ and $\eta$ a natural transformation from $F$ to $F'\circ G^{\mathrm{op}}$, where $G^{\mathrm{op}}$ is the dual of $G$. We let $\mathsf{SICat}$  stand for the corresponding category.

An \emph{opposite split} $\mathsf{I}$-\emph{indexed category} is a covariant functor $F$ from $\mathsf{I}$ to $\mathsf{Cat}$. Given an object $i\in \mathrm{Ob}(\mathsf{I})$, we write $\mathsf{F}_{i}$ for the category $F(i)$, and given a morphism $\varphi\in \mathrm{Hom}_{\mathsf{I}}(i,j)$, we write $F_{\varphi}$ for the functor $F(\varphi)\colon \mathsf{F}_{i}\mor \mathsf{F}_{j}$. An \emph{opposite split indexed category} is an ordered pair $(\mathsf{I},F)$ in which $\mathsf{I}$ is a category and $F$ an opposite split $\mathsf{I}$-indexed category. A \emph{morphism of opposite split indexed categories} from an opposite split indexed category $(\mathsf{I},F)$ to another $(\mathsf{I}',F')$ is an ordered pair $(G,\eta)$, where $G$ is a functor from $\mathsf{I}$ to $\mathsf{I}'$ and $\eta$ a natural transformation from $F$ to $F'\circ G$. We let $\mathsf{OpSICat}$ stand for the corresponding category.
\end{definition}

We next define the object part of the Grothendieck construction (for the functor from the category of split indexed categories to the category of split fibrations) which assigns to a split indexed category $(\mathsf{I},F)$, where $F$ is a split $\mathsf{I}$-indexed category, the split fibration $(\int^{\mathsf{I}}F,\pi_{F})$---and in doing so we follow the notation used in~\cite{JY21} (Chapter~10) but adding the domain of $F$ as a superscript to indicate the contravariance of $F$.

\begin{definition}
Let $(\mathsf{I},F)$ be a split indexed category. Then the \emph{Grothen\-dieck construction} at $(\mathsf{I},F)$ is $(\int^{\mathsf{I}}F,\pi_{F})$, where $\int^{\mathsf{I}}F$ is the category defined as follows:
\begin{enumerate}
\item $\mathrm{Ob}(\int^{\mathsf{I}}F) = \bigcup_{i\in\mathrm{Ob}(\mathsf{I})}(\{i\}\times\mathrm{Ob}(\mathsf{F}_{i}))$.
\item For every $(i,x), (j,y)\in \mathrm{Ob}(\int^{\mathsf{I}}F)$, $\mathrm{Hom}_{\int^{\mathsf{I}}F}((i,x),(j,y))$ is the set of all ordered pairs $(\varphi,f)$, where $\varphi\in \mathrm{Hom}_{\mathsf{I}}(i,j)$ and $f\in \mathrm{Hom}_{\mathsf{F}_{i}}(x,F_{\varphi}(y))$.
\item For every $(i,x)\in \mathrm{Ob}(\int^{\mathsf{I}}F)$, the identity morphism at $(i,x)$ is given by $(\mathrm{id}_{i},\mathrm{id}_{x})$.
\item For every $(i,x), (j,y), (k,z)\in \mathrm{Ob}(\int^{\mathsf{I}}F)$, every $(\varphi,f)\colon (i,x)\mor (j,y)$, and every $(\psi,g)\colon (j,y)\mor (k,z)$, the composite morphism $(\psi,g)\circ (\varphi,f)$ from $(i,x)$ to $(k,z)$ is
     $$
      (\psi,g)\circ (\varphi,f) = (\psi\circ \varphi,F_{\varphi}(g)\circ f)\colon (i,x)\mor (k,z).
     $$
Notice that $f\colon x\mor F_{\varphi}(y)$, $g\colon y\mor F_{\psi}(z)$, hence, taking into account that $F_{\varphi}$ is a functor from $\mathsf{F}_{j}$ to $\mathsf{F}_{i}$, $F_{\varphi}(g)\colon F_{\varphi}(y)\mor F_{\varphi}(F_{\psi}(z))$. Therefore $F_{\varphi}(g)\circ f\colon x\mor F_{\varphi}(F_{\psi}(z))$.
\end{enumerate}
And $\pi_{F}$ the \emph{canonical split fibration} from $\int^{\mathsf{I}}F$ to $\mathsf{I}$ (determined by $(\mathsf{I},F)$) that sends $(\varphi,f)\colon (i,x)\mor (j,y)$ to $\varphi\colon i\mor j$. For details see~\cite{BW12}.
\end{definition}


We next define the morphism part of the Grothendieck construction (for the functor from the category of split indexed categories to the category of split fibrations). But before doing so, we prove that every morphism between split indexed categories determines a morphism between their associated split fibrations.

\begin{proposition}\label{Intmorf}
Let $(G,\eta)$ be a morphism from a split indexed category $(\mathsf{I},F)$ to another  $(\mathsf{I}',F')$. Then there exists a functor $\int^{G}\eta$ from $\int^{\mathsf{I}}F$ to $\int^{\mathsf{I}'}F'$ such that $G\circ \pi_{F} = \pi_{F'}\circ \int^{G}\eta$, i.e., such that the following diagram:
$$\xymatrix{
\int^{\mathsf{I}}F\ar[d]_{\pi_{F}}
\ar[r]^-{\int^{G}\eta} & \int^{\mathsf{I}'}F'\ar[d]^{\pi_{F'}}  \\
\mathsf{I}\ar[r]_-{G} &
\mathsf{I}'
}
$$
commutes.
\end{proposition}

\begin{proof}
If $(i,x)$ is an object of $\int^{\mathsf{I}}F$, then, taking into account that $\eta_{i}$ is a functor from $\mathsf{F}_{i}$ to $\mathsf{F}'_{G(i)}$, we define $(\int^{G}\eta)(i,x)$ as $(G(i),\eta_{i}(x))$. On the other hand, taking into account that the following diagram:
$$\xymatrix{
\mathsf{F}_{i} \ar[r]^-{\eta_{i}} & \mathsf{F}'_{G(i)} \\
\mathsf{F}_{j} \ar[u]^-{F_{\varphi}} \ar[r]_-{\eta_{j}} &
\mathsf{F}'_{G(j)} \ar[u]_-{F'_{G(\varphi)}}
}
$$
commutes, if $(\varphi,f)$ is a morphism from $(i,x)$ to $(j,y)$, then we define the morphism $(\int^{G}\eta)(\varphi,f)$ from $(G(i),\eta_{i}(x))$ to $(G(j),\eta_{j}(y))$ as  $(G(\varphi),\eta_{i}(f))$ (notice that $f\colon x\mor F_{\varphi}(y)$ and $\eta_{i}(F_{\varphi}(y)) = F'_{G(\varphi)}(\eta_{j}(y))$).
\end{proof}

\begin{definition}
Let $(G,\eta)$ be a morphism from a split indexed category $(\mathsf{I},F)$ to another $(\mathsf{I}',F')$. Then the \emph{Grothendieck construction} at $(G,\eta)$ is the morphism $(\int^{G}\eta, G)$ from the split fibration $(\int^{\mathsf{I}}F,\pi_{F})$ to the split fibration $(\int^{\mathsf{I}'}F',\pi_{F'})$. 
\end{definition}

\begin{remark}
There exists a functor $\int^{\bigcdot}$ from the category of split indexed categories  to the category of split fibrations.
\end{remark}

For later use, and for the reader's convenience, we define the object part of the Grothendieck construction (for the functor from the category of opposite split indexed categories to the category of split opfibrations) which assigns to an opposite split indexed category $(\mathsf{I},F)$, where $F$ is an opposite split $\mathsf{I}$-indexed category, the split opfibration $(\int_{\mathsf{I}}F,\pi_{F})$---note that we add the domain of $F$ as a subscript to indicate the covariance of $F$.

\begin{definition}\label{DGConsOp}
Let $(\mathsf{I},F)$ be an opposite split indexed category. Then the \emph{Grothen\-dieck construction} at $(\mathsf{I},F)$ is $(\int_{\mathsf{I}}F,\pi_{F})$, where $\int_{\mathsf{I}}F$ is the category defined as follows:
\begin{enumerate}
\item $\mathrm{Ob}(\int_{\mathsf{I}}F) = \bigcup_{i\in\mathrm{Ob}(\mathsf{I})}(\{i\}\times\mathrm{Ob}(\mathsf{F}_{i}))$.
\item For every $(i,x), (j,y)\in \mathrm{Ob}(\int_{\mathsf{I}}F)$, $\mathrm{Hom}_{\int_{\mathsf{I}}F}((i,x),(j,y))$ is the set of all ordered pairs $(\varphi,f)$, where $\varphi\in \mathrm{Hom}_{\mathsf{I}}(i,j)$ and $f\in \mathrm{Hom}_{\mathsf{F}_{j}}(F_{\varphi}(x),y)$.
\item For every $(i,x)\in \mathrm{Ob}(\int_{\mathsf{I}}F)$, the identity morphism at $(i,x)$ is given by $(\mathrm{id}_{i},\mathrm{id}_{x})$.
\item For every $(i,x), (j,y), (k,z)\in \mathrm{Ob}(\int_{\mathsf{I}}F)$, every $(\varphi,f)\colon (i,x)\mor (j,y)$, and every $(\psi,g)\colon (j,y)\mor (k,z)$, the composite morphism $(\psi,g)\circ (\varphi,f)$ from $(i,x)$ to $(k,z)$ is
     $$
      (\psi,g)\circ (\varphi,f) = (\psi\circ \varphi,g\circ F_{\psi}(f))\colon (i,x)\mor (k,z).
     $$
     Notice that $f\colon F_{\varphi}(x)\mor y$, $g\colon F_{\psi}(y)\mor z$, hence, taking into account that $F_{\psi}$ is a functor from $\mathsf{F}_{j}$ to $\mathsf{F}_{k}$, $F_{\psi}(f)\colon F_{\psi}(F_{\varphi}(x))\mor F_{\psi}(y)$. Therefore $g\circ F_{\psi}(f)\colon F_{\psi}(F_{\varphi}(x))\mor z$.
\end{enumerate}
And $\pi_{F}$ the \emph{canonical split opfibration} from $\int_{\mathsf{I}}F$ to $\mathsf{I}$ (determined by $(\mathsf{I},F)$)  that sends $(\varphi,f)\colon (i,x)\mor (j,y)$ to $\varphi\colon i\mor j$. For details see~\cite{BW12}.
\end{definition}



\begin{proposition}
Let $(G,\eta)$ be a morphism from an opposite split indexed category $(\mathsf{I},F)$ to another $(\mathsf{I}',F')$. Then there exists a functor $\int_{G}\eta$ from $\int_{\mathsf{I}}F$ to $\int_{\mathsf{I}'}F'$ such that $G\circ \pi_{F} = \pi_{F'}\circ \int_{G}\eta$.
\end{proposition}

\begin{remark}
There exists a functor $\int_{\bigcdot}$ from the category of opposite split indexed categories to the category of split opfibrations.
\end{remark}

We next define the notion of functor having the property of being a weak right multiadjoint (or, what is equivalent, having a weak left multiadjoint), which is a generalization of the notion of right multiadjoint defined by Diers, among other places, in~\cite{YD84}. We will prove later on that the Lallement functor is a weak right multiadjoint.

\begin{definition}
Let $G$ be a functor from $\mathsf{A}$ to $\mathsf{X}$. We will say that $G$ \emph{is a weak right multiadjoint} or that $G$ \emph{has a weak left multiadjoint} provided that for every object $X$ in $\mathsf{X}$ there exists a set $I$, a family of objects $(A_{i})_{i\in I}$ in $\mathsf{A}$ and a family of morphisms $(f_{i})_{i\in I}$ in $\prod_{i\in I}\mathrm{Hom}_{\mathsf{X}}(X,G(A_{i}))$ such that, for every object $A$ in $\mathsf{A}$ and every morphism $f$ from $X$ to $G(A)$, there exists a non-empty subset $J$ of $I$ such that, for every $j\in J$, there exists a unique morphism $g_{j}$ from $A_{j}$ to $A$ such that the following diagram
\begin{center}
\begin{tikzpicture}
[AClimentM/.style={|-{To [angle'=45, length=5.75pt, width=4pt, round]}},
ACliment/.style={-{To [angle'=45, length=5.75pt, width=4pt, round]}},scale=1,
AClimentD/.style={double equal sign distance,-implies}, scale=1]
            
\node[] (X) 		at (0,0) [] {$X$};
\node[] (GAj) 	at (4,-0) [] {$G(A_{j})$};
\node[] (GA) 	at (4,-2) [] {$G(A)$};
            
\draw[ACliment]  (X) to node [above]	{$f_{j}$} (GAj);
\draw[ACliment]  (X) to node [below left]	{$f$} (GA);
\draw[ACliment]  (GAj) to node [right]	{$G(g_{j})$} (GA);      
\end{tikzpicture}
\end{center}
commutes, i.e., $f = G(g_{j})\circ f_{j}$. 
\end{definition}

\begin{remark}
For Diers instead of $\exists J\subseteq I$, with $J \neq \emptyset$, one has $\exists ! \,i\in I$. 
\end{remark}

\begin{examples}
Linearly ordered sets are multireflective in the category of strictly ordered sets (we recall that a subcategory $\mathsf{A}$ of a category $\mathsf{X}$ is multireflective if the inclusion functor is a right multiadjoint). The category of fields is multireflective in the category of commutative unitary rings. In~\cite{YD81,YD84} the reader will find a large number of examples, order-theoretical, algebraic and topological, of right multiadjoints.
\end{examples}

\section{Lallement systems of $\Sigma$-algebras}

In this section, for a $\Sigma$-algebra $\mathbf{I}$, we define the notion of Lallement system of $\Sigma$-algebras relative to $\mathbf{I}$, which is taken 
from~\cite[Th\'{e}or\`{e}me 2.19, p.~83]{La67} and \cite[Definition~4.5.2, p.~214]{RS02}, and of morphism between Lallement systems of $\Sigma$-algebras relative to $\mathbf{I}$ in such a way as to obtain a category. Moreover, we define a contravariant functor $\mathrm{Lsys}_{\Sigma}$ from the category $\mathsf{Alg}(\Sigma)$ to the category $\mathsf{Cat}$ and from it we obtain, by means of the Grothendieck construction, the category $\int^{\mathsf{Alg}(\Sigma)}\mathrm{Lsys}_{\Sigma}$, of labeled Lallement systems of $\Sigma$-algebras, and the split fibration $\pi_{\mathrm{Lsys}_{\Sigma}}$ from $\int^{\mathsf{Alg}(\Sigma)}\mathrm{Lsys}_{\Sigma}$ to $\mathsf{Alg}(\Sigma)$.

We begin by defining the notion of Lallement system of $\Sigma$-algebras relative to $\mathbf{I}$ and of morphism between them. 

\begin{definition}\label{DLsys}
Let $\mathbf{I}=(I,F^{\mathbf{I}})$ be a $\Sigma$-algebra with algebraic predecessor preorder $\leq_{\mathbf{I}}$. Then a \emph{Lallement system of $\Sigma$-algebras relative to $\mathbf{I}$}, or an \emph{$\mathbf{I}$-Lallement system of $\Sigma$-algebras} for brevity, is an ordered pair 
$\mathcal{A}=((\mathbf{A}_{i}, \mathbf{M}_{i})_{i\in I}, (f_{i,j})_{(i,j)\in\leq_{\mathbf{I}}})$, where $(\mathbf{A}_{i},\mathbf{M}_{i})_{i\in I}$ is an $I$-indexed family of pairs of $\Sigma$-algebras such that, for every $i\in I$, $\mathbf{A}_{i}=(A_{i},F^{\mathbf{A}_{i}})$, $\mathbf{M}_{i}=(M_{i},F^{\mathbf{M}_{i}})$, $\mathbf{A}_{i}\in\mathrm{Sub}(\mathbf{M}_{i})$ and $(f_{i,j})_{(i,j)\in\leq_{\mathbf{I}}}$ a family of homomorphisms in $\prod_{(i,j)\in\leq_{\mathbf{I}}}\mathrm{Hom}(\mathbf{A}_{i},\mathbf{M}_{j})$ such that
\begin{enumerate}
\item[(a)] for every $n\in\mathbb{N}$, every $\sigma\in\Sigma_{n}$, every $(i_{k})_{k\in n}\in I^{n}$ with $F_{\sigma}^{\mathbf{I}}((i_{k})_{k\in n})=i$ and every $(a_{i_{k}})_{k\in n}\in \prod_{k\in n} A_{i_{k}}$,
\[
F_{\sigma}^{\mathbf{M}_{i}}
((f_{i_{k},i}(a_{i_{k}}))_{k\in n})\in A_{i}, \text{ and}
\]
\item[(b)] for every $n\in\mathbb{N}$, every $\sigma\in\Sigma_{n}$, every $(i_{k})_{k\in n}\in I^{n}$ with $F_{\sigma}^{\mathbf{I}}((i_{k})_{k\in n})=i$, every $(a_{i_{k}})_{k\in n}\in \prod_{k\in n} A_{i_{k}}$ and every $j\in I$, if $(i,j)\in\leq_{\mathbf{I}}$, then
\[
f_{i,j}(F_{\sigma}^{\mathbf{M}_{i}}((f_{i_{k},i}(a_{i_{k}}))_{k\in n}))=
F_{\sigma}^{\mathbf{M}_{j}}((f_{i_{k},j}(a_{i_{k}}))_{k\in n}).
\]
\end{enumerate}
We will refer to the ordered pairs $(\mathbf{A}_{i},\mathbf{M}_{i})$ of an $\mathbf{I}$-Lallement system of $\Sigma$-algebras $\mathcal{A}$ as the \emph{algebraic pairs of} $\mathcal{A}$, and to the homomorphisms $f_{i,j}$ as the \emph{transition homomorphisms of} $\mathcal{A}$. 
Let $\mathcal{A}$ and $\mathcal{B}=((\mathbf{B}_{i}, \mathbf{N}_{i})_{i\in I}, (g_{i,j})_{(i,j)\in\leq_{\mathbf{I}}})$ be $\mathbf{I}$-Lallement systems of $\Sigma$-algebras. Then a \emph{morphism} from $\mathcal{A}$ to $\mathcal{B}$ is an ordered triple $(\mathcal{A},u,\mathcal{B})$, abbreviated to $u\colon\mathcal{A}\mor \mathcal{B}$, where $u=(u_{i})_{i\in I}$ is a family of homomorphisms in $\prod_{i\in I}\mathrm{Hom}(\mathbf{M}_{i},\mathbf{N}_{i})$ such that, for every $i\in I$, $u_{i}[A_{i}]\subseteq B_{i}$ and, for every $(i,j)\in\leq_{\mathbf{I}}$, $u_{j}\circ f_{i,j}=g_{i,j}\circ u_{i}$, where, by abuse of notation, we continue to write $u_{i}$ for the birestriction of the homomorphism $u_{i}$ from $\mathbf{M}_{i}$ to $\mathbf{N}_{i}$ to $\mathbf{A}_{i}$ and $\mathbf{B}_{i}$. We denote by $\mathbf{I}$-$\mathsf{LAlg}(\Sigma)$ the corresponding category (which is an object of $\mathsf{Cat}$).
\end{definition}

\begin{remark}
The conditions $(a)$ and $(b)$ in Definition~\ref{DLsys} are the conditions $1)$ and $2)$ stated by Lallement in~\cite[Th\'{e}or\`{e}me 2.19, pp.~83--84]{La67}, as generalized by Romanowska and Smith in~\cite{RS91} to the field of universal algebra. Moreover, $(a)$ can be restated as: 
$$
\textstyle
F_{\sigma}^{\mathbf{M}_{i}}\circ \prod_{k\in n}f_{i_{k},i}\left[\prod_{k\in n}A_{i_{k}}\right]\subseteq A_{i},
$$ 
where $\prod_{k\in n}f_{i_{k},i}$ is the unique mapping from $\prod_{k\in n}A_{i_{k}}$ to $M^{n}_{i}$ such that, for every $k\in n$, 
$f_{i_{k},i}\circ \mathrm{pr}_{i_{k}} = \mathrm{pr}_{k}\circ \prod_{k\in n}f_{i_{k},i}$, with $\mathrm{pr}_{i_{k}}$ the canonical projection from $\prod_{k\in n}A_{i_{k}}$ to $A_{i_{k}}$ and $\mathrm{pr}_{k}$ the canonical projection from $M^{n}_{i}$ to $M_{i}$, for every $k\in n$; and $(b)$ can be restated as: 
$$
\textstyle
f_{i,j}\circ\mathrm{cor}(F_{\sigma}^{\mathbf{M}_{i}}\circ \prod_{k\in n}f_{i_{k},i})_{A_{i}} = F_{\sigma}^{\mathbf{M}_{j}}\circ \prod_{k\in n}f_{i_{k},j},
$$ 
where $\mathrm{cor}(F_{\sigma}^{\mathbf{M}_{i}}\circ \prod_{k\in n}f_{i_{k},i})_{A_{i}}$ is the corestriction of $F_{\sigma}^{\mathbf{M}_{i}}\circ \prod_{k\in n}f_{i_{k},i}$ to $A_{i}$ (which exists by $(a)$).
\end{remark}

We next define the assignment $\mathrm{Lsys}_{\Sigma}$ from $\mathsf{Alg}(\Sigma)$ to $\mathsf{Cat}$. 

\begin{definition}\label{DLsysFun} 
We let $\mathrm{Lsys}_{\Sigma}$ stand for the assignment from $\mathsf{Alg}(\Sigma)$ to $\mathsf{Cat}$ defined as follows:
\begin{enumerate}
\item for every $\Sigma$-algebra $\mathbf{I}$, $\mathrm{Lsys}_{\Sigma}(\mathbf{I})$ is $\mathbf{I}$-$\mathsf{LAlg}(\Sigma)$, and
\item for every homomorphism $t\colon \mathbf{I}\mor \mathbf{P}$, $\mathrm{Lsys}_{\Sigma}(t)$ is the functor from $\mathbf{P}$-$\mathsf{LAlg}(\Sigma)$ to $\mathbf{I}$-$\mathsf{LAlg}(\Sigma)$ defined as follows:
\begin{enumerate}
\item for every $\mathbf{P}$-Lallement system of $\Sigma$-algebras $\mathcal{A}=((\mathbf{A}_{p}, \mathbf{M}_{p})_{p\in P}, (f_{p,q})_{(p,q)\in \leq_{\mathbf{P}}})$, $\mathrm{Lsys}_{\Sigma}(t)(\mathcal{A})$, denoted by $\mathcal{A}_{t}$ for short, is $((\mathbf{A}_{t(i)}, \mathbf{M}_{t(i)})_{i\in I}, (f_{t(i),t(j)})_{(i,j)\in \leq_{\mathbf{I}}})$, which is an $\mathbf{I}$-Lallement system of $\Sigma$-algebras, and
\item for every morphism of $\mathbf{P}$-Lallement systems of $\Sigma$-algebras $u$, where $u=(u_{p})_{p\in P}$, from $\mathcal{A}$ to $\mathcal{B}=((\mathbf{B}_{p}, \mathbf{N}_{p})_{p\in P}, (g_{p,q})_{(p,q)\in \leq_{\mathbf{P}}})$, $\mathrm{Lsys}_{\Sigma}(t)(u)$, denoted by $u_{t}$ for short, is $(u_{t(i)})_{i\in I}$, which is a morphism of $\mathbf{I}$-Lallement system of $\Sigma$-algebras from $\mathcal{A}_{t}$ to $\mathcal{B}_{t}$.
\end{enumerate}
\end{enumerate}
\end{definition}

\begin{proposition} 
The assignment $\mathrm{Lsys}_{\Sigma}$ from $\mathsf{Alg}(\Sigma)$ to $\mathsf{Cat}$ is a contravariant functor.
\end{proposition}

\begin{definition}\label{DLsysGC} 
From the split indexed category $(\mathsf{Alg}(\Sigma), \mathrm{Lsys}_{\Sigma})$ we obtain, by means of the Grothendieck construction, 
\begin{enumerate}
\item the category $\int^{\mathsf{Alg}(\Sigma)}\mathrm{Lsys}_{\Sigma}$ which has as objects the \emph{Lallement systems of $\Sigma$-algebras labeled by $\Sigma$-algebras}, i.e., the ordered pairs $\mathbcal{A}=(\mathbf{I},\mathcal{A})$, where $\mathbf{I}$ is a $\Sigma$-algebra and $\mathbcal{A}$ an $\mathbf{I}$-Lallement system of $\Sigma$-algebras; and as morphisms from the object $\mathbcal{A}$ to the object $\mathbcal{B}=(\mathbf{P},\mathcal{B})$ the ordered triples $(\mathbcal{A}, (t,u), \mathbcal{B})$, abbreviated to $(t,u)\colon \mathbcal{A}\mor \mathbcal{B}$, where $t$ is a homomorphism from $\mathbf{I}$ to $\mathbf{P}$ and $u$ a morphism from the $\mathbf{I}$-Lallement system of $\Sigma$-algebras $\mathcal{A}$ to the $\mathbf{I}$-Lallement system of $\Sigma$-algebras $\mathcal{B}_{t}$; and 
\item the split fibration $\pi_{\mathrm{Lsys}_{\Sigma}}$ from $\int^{\mathsf{Alg}(\Sigma)}\mathrm{Lsys}_{\Sigma}$ to $\mathsf{Alg}(\Sigma)$ that sends $\mathbcal{A}=(\mathbf{I},\mathcal{A})$ to $\mathbf{I}$ and $(t,u)\colon \mathbcal{A}\mor \mathbcal{B}$, with $\mathbcal{B}=(\mathbf{P},\mathcal{B})$ to $t\colon \mathbf{I}\mor \mathbf{P}$.
\end{enumerate}
\end{definition}

\section{The category of naturally preordered idempotent $\Sigma$-algebras and its relatives}

In this section, after defining the notion of plural signature, we define the notion of idempotent $\Sigma$-algebra and provide a characterization of the idempotency of a $\Sigma$-algebra in terms of constant mappings. After that, we define the concept of naturally preordered $\Sigma$-algebra, the category $\mathsf{NPIAlg}(\Sigma)$, of naturally preordered idempotent $\Sigma$-algebras, and provide several characterizations of the concept of naturally preordered idempotent $\Sigma$-algebra. 
Moreover, from the category $\mathsf{NPIAlg}(\Sigma)_{\mathsf{s}}$, of naturally preordered idempotent $\Sigma$-algebras and surjective homomorphisms, we obtain, on the one hand, by restricting
the contravariant functor $\mathrm{Lsys}_{\Sigma}$ to it, the category 
$\int^{\mathsf{NPIAlg}(\Sigma)_{\mathsf{s}}}\mathrm{Lsys}_{\Sigma}$ which has as objects the Lallement systems of $\Sigma$-algebras labeled by naturally preordered idempotent $\Sigma$-algebras; and, on the other hand, after defining another functor 
$(\mathsf{Alg}(\Sigma)\,{\downarrow_{\mathsf{s}}}\, \bigcdot)$ from $\mathsf{NPIAlg}(\Sigma)_{\mathsf{s}}$ to $\mathsf{Cat}$, the category $\int_{\mathsf{NPIAlg}(\Sigma)_{\mathsf{s}}}(\mathsf{Alg}(\Sigma)\,{\downarrow_{\mathsf{s}}}\, \bigcdot)$ which has as objects the ordered pairs 
$(\mathbf{I},(\mathbf{A},f))$, where $\mathbf{I}$ is a naturally preordered idempotent $\Sigma$-algebra and $(\mathbf{A},f)$ an object of $(\mathsf{Alg}(\Sigma)\,{\downarrow_{\mathsf{s}}}\, \mathbf{I})$, i.e., a $\Sigma$-algebra 
$\mathbf{A}$ and a surjective homomorphism $f$ from $\mathbf{A}$ to $\mathbf{I}$. We point out that 
the Lallement functor that will have the property of being a weak right multiadjoint will be the one that has as domain a suitable subcategory of $\int^{\mathsf{NPIAlg}(\Sigma)_{\mathsf{s}}}\mathrm{Lsys}_{\Sigma}$ and as codomain the category $\int_{\mathsf{NPIAlg}(\Sigma)_{\mathsf{s}}}(\mathsf{Alg}(\Sigma)\,{\downarrow_{\mathsf{s}}}\, \bigcdot)$.

\begin{definition}\label{DPSig}
A signature $\Sigma$ will be called \emph{plural} if $\Sigma_{0}=\varnothing$ and, for some $n\in\mathbb{N}$ with $n\geq 2$, $\Sigma_{n}\neq\varnothing$.
\end{definition}


\begin{assumption}\label{ALal}
From now on, unless otherwise stated, $\Sigma$ stands for a plural signature as set forth in Definition~\ref{DPSig}, fixed once and for all. In addition, no $\Sigma$-algebra will be empty.
\end{assumption}

We begin by defining the notion of idempotent $\Sigma$-algebra.

\begin{definition}\label{DIdem} 
Let $\mathbf{I}=(I,F)$ be a $\Sigma$-algebra, $n\in \mathbb{N}-\{0\}$ and $\sigma\in \Sigma_{n}$. We will say that $F^{\mathbf{I}}_{\sigma}$ is \emph{idempotent} if, for every $i\in I$, $F^{\mathbf{I}}_{\sigma}((i)_{k\in n})=i$, where $(i)_{k\in n}$ is the family in $I^{n}$ which is constantly $i$. We will say that $\mathbf{I}$ is \emph{idempotent} if, for every $n\in\mathbb{N}-\{0\}$ and every $\sigma\in \Sigma_{n}$, $F^{\mathbf{I}}_{\sigma}$ is idempotent. We will denote by $\mathsf{IAlg}(\Sigma)$ the full subcategory of $\mathsf{Alg}(\Sigma)$ determined by the idempotent $\Sigma$-algebras.
\end{definition}

We next provide a characterization of the idempotent $\Sigma$-algebras by means of the constant mappings.

\begin{proposition}\label{PIdem} 
Let $\mathbf{I}=(I,F)$ be a $\Sigma$-algebra. Then the following statements are equivalent.
\begin{enumerate}
\item $\mathbf{I}$ is idempotent.
\item Every constant mapping $\kappa\colon I\mor I$ is an endomorphism.
\end{enumerate}
\end{proposition}

%
%

\begin{remark}
Let $\mathbf{I}=(I,F)$ and $\mathbf{P}=(P,G)$ be idempotent $\Sigma$-algebras. Then, for every $p\in P$, the mapping $\kappa_{p}$ from $I$ to $P$ which is constantly $p$, determines a homomorphism from $\mathbf{I}$ to $\mathbf{P}$. 
\end{remark}

We next define the concept of naturally preordered $\Sigma$-algebra.

\begin{definition}\label{DNatPOrd} 
Let $\mathbf{I}=(I,F)$ be a $\Sigma$-algebra. We say that $\mathbf{I}$ is \emph{naturally preordered} if the algebraic predecessor preorder $\leq_{\mathbf{I}}$ on $\mathbf{I}$ is such that, for every $n\in\mathbb{N}-\{0\}$ and every $\sigma\in\Sigma_{n}$, if $(j_{k})_{k\in n}$ and $(i_{k})_{k\in n}$ are families in $I^{n}$ such that, for every $k\in n$, $j_{k}\leq_{\mathbf{I}} i_{k}$, then $F^{\mathbf{I}}_{\sigma}((j_{k})_{k\in n})\leq_{\mathbf{I}} F^{\mathbf{I}}_{\sigma}((i_{k})_{k\in n})$.
\end{definition}

\begin{remark}
A $\Sigma$-algebra $\mathbf{I}$ is naturally preordered if, and only if, for every $n\in\mathbb{N}-\{0\}$ and every $\sigma\in\Sigma_{n}$, $F^{\mathbf{I}}_{\sigma}$ is an isotone mapping from 
$(I^{n},\leq_{\mathbf{I}^{n}})$ to $(I,\leq_{\mathbf{I}})$, where $\leq_{\mathbf{I}^{n}}$ is the preorder on $I^{n}$ defined, for every $(j_{k})_{k\in n}$, $(i_{k})_{k\in n}\in I^{n}$ as follows: 
$$
(j_{k})_{k\in n}\leq_{\mathbf{I}^{n}}(i_{k})_{k\in n} \text{ if, and only if, for every } k\in n,\, 
j_{k}\leq_{\mathbf{I}}i_{k}.
$$
\end{remark}

\begin{proposition}
Let $f$ be a homomorphism from the naturally preordered $\Sigma$-algebra $\mathbf{I}=(I,F)$ to the 
naturally preordered $\Sigma$-algebra $\mathbf{P}=(P,G)$. Then, for every $n\in\mathbb{N}-\{0\}$ and every $\sigma\in\Sigma_{n}$, if $(j_{k})_{k\in n}$ and $(i_{k})_{k\in n}$ are families in $I^{n}$ such that, for every $k\in n$, $j_{k}\leq_{\mathbf{I}} i_{k}$, then 
$f(F_{\sigma}((j_{k})_{k\in n}))\leq_{\mathbf{P}} f(F_{\sigma}((i_{k})_{k\in n}))$.
\end{proposition}

We next introduce the concept of $\Sigma$-semilattice, which will be used afterwards when we provide characterizations of the notion of naturally preordered idempotent $\Sigma$-algebra. 

\begin{definition}\label{DSSsl} 
Let $\mathbf{I}=(I,F)$ be a $\Sigma$-algebra. We will say that $\mathbf{I}$ is a \emph{$\Sigma$-semilattice} if there exists an order $\leq$ on $I$ such that $(I,\leq)$ is a semilattice and, for every $n\in \mathbb{N}-\{0\}$, every $\sigma\in\Sigma_{n}$ and every $(i_{k})_{k\in n}\in I^{n}$, we have that $\bigvee_{k\in n}i_{k}=F_{\sigma}((i_{k})_{k\in n})$.
\end{definition}

In the following proposition we provide several characterizations of the concept of naturally preordered idempotent $\Sigma$-algebra.

\begin{proposition}\cite[Prop.~4.1.7, pp.~184--185]{RS02}\label{PIdemNatPOrd} 
Let $\mathbf{I}=(I,F)$ be an idempotent $\Sigma$-algebra with algebraic predecessor preorder $\leq_{\mathbf{I}}$. Then the following statements are equivalent:
\begin{enumerate}
\item $\mathbf{I}$ is naturally preordered.
\item For every $n\in\mathbb{N}-\{0\}$, every $\sigma \in \Sigma_{n}$, every $i\in I$ and every $(i_{k})_{k\in n}$ in $I^{n}$, if, for every $k\in n$, $i_{k}\leq_{\mathbf{I}} i$, then $F_{\sigma}((i_{k})_{k\in n})\leq_{\mathbf{I}} i$.
\item The relation $\equiv_{\mathbf{I}}=\leq_{\mathbf{I}}\cap \geq_{\mathbf{I}}$ is a congruence on $\mathbf{I}$ and the quotient $\mathbf{I}/{\equiv_{\mathbf{I}}}$ is a $\Sigma$-semilattice.
\end{enumerate}
\end{proposition}

\begin{convention}
From now on, for brevity, instead of saying ``naturally preordered idempotent $\Sigma$-algebra'' we will say ``npi $\Sigma$-algebra''. 
\end{convention}

\begin{corollary}\label{CIdemNat} 
Let $\mathbf{I}=(I,F)$ be an npi $\Sigma$-algebra. Then, for every $i\in I$, $[i]_{\equiv_{\mathbf{I}}}$ is 
a closed subset of $\mathbf{I}$. In particular, $\mathbf{I}$ is a union of subalgebras.
\end{corollary}

%
%
%
%

\begin{definition}\label{DNOIdpAlg} 
We will denote by $\mathsf{NPIAlg}(\Sigma)$ the full subcategory of $\mathsf{IAlg}(\Sigma)$ determined by the naturally preordered $\Sigma$-algebras. Moreover, we will denote by $\mathsf{NPIAlg}(\Sigma)_{\mathsf{s}}$ the subcategory of $\mathsf{NPIAlg}(\Sigma)$ consisting of the npi $\Sigma$-algebras and the surjective homomorphisms between them.
\end{definition}


\begin{definition}\label{DLsysGC2}
By Definition~\ref{DNOIdpAlg}, we can consider the restriction of the contravariant functor $\mathrm{Lsys}_{\Sigma}$, set forth in Definition~\ref{DLsysFun}, to $\mathsf{NPIAlg}(\Sigma)_{\mathsf{s}}$ which, for simplicity's sake, we will also denote by $\mathrm{Lsys}_{\Sigma}$. From the split indexed category $(\mathsf{NPIAlg}(\Sigma)_{\mathsf{s}}, \mathrm{Lsys}_{\Sigma})$ we obtain, by means of the Grothendieck construction, 
\begin{enumerate}
\item the category $\int^{\mathsf{NPIAlg}(\Sigma)_{\mathsf{s}}}\mathrm{Lsys}_{\Sigma}$ which has as objects the \emph{Lallement systems of $\Sigma$-algebras labeled by npi $\Sigma$-algebras}, i.e., the ordered pairs $\mathbcal{A}=(\mathbf{I},\mathcal{A})$, where $\mathbf{I}$ is an npi $\Sigma$-algebra and $\mathbcal{A}$ an $\mathbf{I}$-Lallement system of $\Sigma$-algebras; and as morphisms from the object $\mathbcal{A}$ to the object $\mathbcal{B}=(\mathbf{P},\mathcal{B})$ the ordered triples $(\mathbcal{A}, (t,u), \mathbcal{B})$, abbreviated to $(t,u)\colon \mathbcal{A}\mor \mathbcal{B}$, where $t$ is a surjective homomorphism from $\mathbf{I}$ to $\mathbf{P}$ and $u$ a morphism from the $\mathbf{I}$-Lallement system of $\Sigma$-algebras $\mathcal{A}$ to the $\mathbf{I}$-Lallement system of $\Sigma$-algebras $\mathcal{B}_{t}$; and 
\item the split fibration $\pi_{\mathrm{Lsys}_{\Sigma}}$ from $\int^{\mathsf{NPIAlg}(\Sigma)_{\mathsf{s}}}\mathrm{Lsys}_{\Sigma}$ to $\mathsf{NPIAlg}(\Sigma)_{\mathsf{s}}$ that sends $\mathbcal{A}=(\mathbf{I},\mathcal{A})$ to $\mathbf{I}$ and $(t,u)\colon \mathbcal{A}\mor \mathbcal{B}$, with $\mathbcal{B}=(\mathbf{P},\mathcal{B})$ to $t\colon \mathbf{I}\mor \mathbf{P}$. 
\end{enumerate}
\end{definition}

We next define the assignment $(\mathsf{Alg}(\Sigma)\,{\downarrow_{\mathsf{s}}}\, \bigcdot)$ from $\mathsf{NPIAlg}(\Sigma)_{\mathsf{s}}$ to $\mathsf{Cat}$.

\begin{definition}\label{DAlgda} 
We let $(\mathsf{Alg}(\Sigma)\,{\downarrow_{\mathsf{s}}}\, \bigcdot)$ stand for the assignment from $\mathsf{NPIAlg}(\Sigma)_{\mathsf{s}}$ to $\mathsf{Cat}$ defined as follows:
\begin{enumerate}
\item for every npi $\Sigma$-algebra $\mathbf{I}$, $(\mathsf{Alg}(\Sigma)\,{\downarrow_{\mathsf{s}}}\, \bigcdot)(\mathbf{I})$, denoted by $(\mathsf{Alg}(\Sigma)\,{\downarrow_{\mathsf{s}}}\, \mathbf{I})$ for short, is the category whose objects are the coverings of $\mathbf{I}$, i.e., the ordered pairs 
$(\mathbf{A},f)$ in which $\mathbf{A}$ is a $\Sigma$-algebra and $f\colon \mathbf{A}\mor \mathbf{I}$ a surjective homomorphism; and as morphisms from $(\mathbf{A},f)$ to $(\mathbf{B},g)$ the homomorphisms $h\colon \mathbf{A}\mor \mathbf{B}$ such that $f=g\circ h$, and
\item for every surjective homomorphism $t\colon \mathbf{I}\mor \mathbf{P}$, $(\mathsf{Alg}(\Sigma)\,{\downarrow_{\mathsf{s}}}\, \bigcdot)(t)$, denoted by $(\mathsf{Alg}(\Sigma)\,{\downarrow_{\mathsf{s}}}\, t)$ for short, is the functor from $(\mathsf{Alg}(\Sigma)\,{\downarrow_{\mathsf{s}}}\, \mathbf{I})$ to $(\mathsf{Alg}(\Sigma)\,{\downarrow_{\mathsf{s}}}\, \mathbf{P})$ defined as follows:
\begin{enumerate}
\item for every object $(\mathbf{A},f)$ in $(\mathsf{Alg}(\Sigma)\,{\downarrow_{\mathsf{s}}}\, \mathbf{I})$, $(\mathsf{Alg}(\Sigma)\,{\downarrow_{\mathsf{s}}}\, t)(\mathbf{A},f)$ is given by $(\mathbf{A}, t\circ f)$, which is an object in $(\mathsf{Alg}(\Sigma)\,{\downarrow_{\mathsf{s}}}\, \mathbf{P})$, and 
\item for every morphism $h\colon (\mathbf{A},f)\mor (\mathbf{B},g)$ in $(\mathsf{Alg}(\Sigma)\,{\downarrow_{\mathsf{s}}}\, \mathbf{I})$, i.e., a homomorphism $h\colon \mathbf{A}\mor \mathbf{B}$ such that $f=g\circ h$, $(\mathsf{Alg}(\Sigma)\,{\downarrow_{\mathsf{s}}}\, t)(h)$ is given by $h$ itself, which is a morphism in $(\mathsf{Alg}(\Sigma)\,{\downarrow_{\mathsf{s}}}\, \mathbf{P})$ from $(\mathbf{A}, t\circ f)$ to $(\mathbf{B}, t\circ g)$.
\end{enumerate}
\end{enumerate}
\end{definition}

\begin{proposition} 
The assignment $(\mathsf{Alg}(\Sigma)\,{\downarrow_{\mathsf{s}}}\, \bigcdot)$ from $\mathsf{NPIAlg}(\Sigma)_{\mathsf{s}}$ to $\mathsf{Cat}$ is a covariant functor.
\end{proposition}

\begin{definition}\label{DAlgdGC} 
From the opposite split indexed category $(\mathsf{NPIAlg}(\Sigma)_{\mathsf{s}}, (\mathsf{Alg}(\Sigma)\,{\downarrow_{\mathsf{s}}}\, \bigcdot))$ we obtain, by means of the Grothendieck construction, 
\begin{enumerate}
\item the category $\int_{\mathsf{NPIAlg}(\Sigma)_{\mathsf{s}}}(\mathsf{Alg}(\Sigma)\,{\downarrow_{\mathsf{s}}}\, \bigcdot)$ which has as objects the ordered pairs $(\mathbf{I},(\mathbf{A},f))$, where $\mathbf{I}$ is an npi $\Sigma$-algebra and $(\mathbf{A},f)$ an object of $(\mathsf{Alg}(\Sigma)\,{\downarrow_{\mathsf{s}}}\, \mathbf{I})$; and as morphisms from the object $(\mathbf{I},(\mathbf{A},f))$ to the object $(\mathbf{P},(\mathbf{B},g))$ the ordered triples $((\mathbf{I},(\mathbf{A},f)), (t, h), (\mathbf{P},(\mathbf{B},g)))$, abbreviated to 
$$
(t,h)\colon (\mathbf{I},(\mathbf{A},f))\mor (\mathbf{P},(\mathbf{B},g)) \text{ for short},
$$ where $t$ is a surjective homomorphism from $\mathbf{I}$ to $\mathbf{P}$ and $h\colon (\mathbf{A},t\circ f)\mor (\mathbf{B}, g)$ a morphism in $(\mathsf{Alg}(\Sigma)\,{\downarrow_{\mathsf{s}}}\, \mathbf{P})$, i.e., a homomorphism $h\colon \mathbf{A}\mor \mathbf{B}$ such that $t\circ f=g\circ h$; and
\item the opposite split fibration 
$\pi_{(\mathsf{Alg}(\Sigma)\,{\downarrow_{\mathsf{s}}}\, \bigcdot)}$ from $\int_{\mathsf{NPIAlg}(\Sigma)_{\mathsf{s}}}(\mathsf{Alg}(\Sigma)\,{\downarrow_{\mathsf{s}}}\, \bigcdot)$ to 
$\mathsf{NPIAlg}(\Sigma)_{\mathsf{s}}$ that sends $(\mathbf{I},(\mathbf{A},f))$ to $\mathbf{I}$ and $(t,h)\colon (\mathbf{I},(\mathbf{A},f))\mor (\mathbf{P},(\mathbf{B},g))$ to $t\colon \mathbf{I}\mor \mathbf{P}$. 
\end{enumerate}
\end{definition}

\section{The functor of Lallement}

In this section, for a plural signature $\Sigma$, we define a functor $\mathfrak{L}_{\Sigma}$, the functor of Lallement, from the category 
$\int^{\mathsf{NPIAlg}(\Sigma)_{\mathsf{s}}}\mathrm{Lsys}_{\Sigma}$ to the category
$\int_{\mathsf{NPIAlg}(\Sigma)_{\mathsf{s}}}(\mathsf{Alg}(\Sigma)\,{\downarrow_{\mathsf{s}}}\, \bigcdot)$.

We begin by defining the notion of Lallement sum of a Lallement system of $\Sigma$-algebras relative to an  npi $\Sigma$-algebra.

\begin{definition}\label{DLaSum}
Let $\mathbcal{A}=(\mathbf{I},\mathcal{A})$ be an object of $\int^{\mathsf{NPIAlg}(\Sigma)_{\mathsf{s}}}\mathrm{Lsys}_{\Sigma}$, where $\mathbf{I}$ is an npi $\Sigma$-algebra and $\mathcal{A}=((\mathbf{A}_{i},\mathbf{M}_{i})_{i\in I}, (f_{i,j})_{(i,j)\in\leq_{\mathbf{I}}})$ an $\mathbf{I}$-Lallement system of $\Sigma$-algebras. The \emph{Lallement sum} of $\mathbcal{A}$ is the $\Sigma$-algebra $\boldsymbol{\mathfrak{L}}_{\Sigma}(\mathbcal{A})$ defined as follows: the underlying set of $\boldsymbol{\mathfrak{L}}_{\Sigma}(\mathbcal{A})$ is $\coprod_{i\in I} A_{i}=\bigcup_{i\in I}(A_{i}\times \{i\})$, the coproduct of $(A_{i})_{i\in I}$, and, for every $n\in\mathbb{N}-\{0\}$, every $\sigma\in\Sigma_{n}$ and every $(a_{i_{k}},i_{k})_{k\in n}\in (\coprod_{i\in I}A_{i})^{n}$,
\[
F^{\boldsymbol{\mathfrak{L}}_{\Sigma}(\mathbcal{A})}_{\sigma}((
a_{i_{k}},i_{k}
)_{k\in n})
=
(F^{\mathbf{M}_{i}}_{\sigma}((f_{i_{k},i}(a_{i_{k}}))_{k\in n}),i),
\]
where $i=F^{\mathbf{I}}_{\sigma}((i_{k})_{k\in n})$. Let us recall that, by condition~(a) of  Definition~\ref{DLsys}, $\boldsymbol{\mathfrak{L}}_{\Sigma}(\mathbcal{A})$ is a well-defined $\Sigma$-algebra.
\end{definition}

\begin{proposition}\label{PLaSum} 
Let $\mathbcal{A}=(\mathbf{I},\mathcal{A})$ be an object of $\int^{\mathsf{NPIAlg}(\Sigma)_{\mathsf{s}}}\mathrm{Lsys}_{\Sigma}$, where $\mathbf{I}$ is an npi $\Sigma$-algebra and $\mathcal{A}$ an $\mathbf{I}$-Lallement system of $\Sigma$-algebras. Then the canonical projection
\[
\mathrm{p}_{\mathbcal{A}}\colon
\left\lbrace
\begin{array}{rcl}
\boldsymbol{\mathfrak{L}}_{\Sigma}(\mathbcal{A})&\mor&\mathbf{I}\\
(a,i)&\longmapsto&i
\end{array}
\right.
\]
is a surjective homomorphism. 
\end{proposition}
\begin{proof}
That $\mathrm{p}_{\mathbcal{A}}$ is surjective follows from the fact that we are assuming that no $\Sigma$-algebra is empty. That $\mathrm{p}_{\mathbcal{A}}$ is an homomorphism follows  from Definition~\ref{DLaSum}.
\end{proof}

%

\begin{proposition}\label{PLaSumMor} 
Let $\mathbcal{A}=(\mathbf{I},\mathcal{A})$ and $\mathbcal{B}=(\mathbf{P},\mathcal{B})$ be two objects of $\int^{\mathsf{NPIAlg}(\Sigma)_{\mathsf{s}}}\mathrm{Lsys}_{\Sigma}$ and $(t,u)$ a morphism in $\int^{\mathsf{NPIAlg}(\Sigma)_{\mathsf{s}}}\mathrm{Lsys}_{\Sigma}$ from $\mathbcal{A}$ to $\mathbcal{B}$. Then there exists a unique homomorphism $u^{t}$ from $\boldsymbol{\mathfrak{L}}_{\Sigma}(\mathbcal{A})$ to $\boldsymbol{\mathfrak{L}}_{\Sigma}(\mathbcal{B})$ 
such that $t\circ \mathrm{p}_{\mathbcal{A}}=\mathrm{p}_{\mathbcal{B}}\circ u^{t}$.
\end{proposition}
\begin{proof}
If $\mathcal{A}=((\mathbf{A}_{i}, \mathbf{M}_{i})_{i\in I}, (f_{i,j})_{(i,j)\in \leq_{\mathbf{I}}})$ and $\mathcal{B}=((\mathbf{B}_{p}, \mathbf{N}_{p})_{p\in P}, (g_{p,q})_{(p,q)\in \leq_{\mathbf{P}}})$, then $t\colon\mathbf{I}\mor \mathbf{P}$ is a surjective homomorphism and  $u\colon\mathcal{A}\mor \mathcal{B}_{t}$ is a morphism of $\mathbf{I}$-Lallement systems of $\Sigma$-algebras where, by  Definition~\ref{DLsysFun}, $\mathcal{B}_{t}$ is the $\mathbf{I}$-Lallement system of $\Sigma$-algebras 
$\mathcal{B}_{t}=((\mathbf{B}_{t(i)}, \mathbf{N}_{t(i)})_{i\in I}, (g_{t(i),t(j)})_{(i,j)\in \leq_{\mathbf{I}}})$. Thus, by Definition~\ref{DLsys}, $u=(u_{i})_{i\in I}$ is a family of homomorphisms in $\prod_{i\in I}\mathrm{Hom}(\mathbf{M}_{i},\mathbf{N}_{t(i)})$ such that, for every $i\in I$, $u_{i}[A_{i}]\subseteq B_{t(i)}$ and, for every $(i,j)\in\leq_{\mathbf{I}}$,
$u_{j}\circ f_{i,j}=g_{t(i),t(j)}\circ u_{i}$. 
For every $i\in I$, let $u^{t}_{i}=\mathrm{in}_{t(i)}\circ u_{i}$ be the mapping from $A_{i}$ to $\coprod_{p\in P}B_{p}$ obtained by composing the birestriction of $u_{i}$ to $A_{i}$ and $B_{t(i)}$, denoted also by $u_{i}$ for short, and $\mathrm{in}_{t(i)}$, the canonical inclusion of $B_{t(i)}$ into $\coprod_{p\in P} B_{p}$. Then, by the universal property of the coproduct, there exists a unique mapping $[(u^{t}_{i})_{i\in I}]\colon \coprod_{i\in I}A_{i}\mor \coprod_{p\in P}B_{p}$, denoted by $u^{t}$ for short, such that, for every $i\in I$, $u^{t}\circ \mathrm{in}_{i}=u^{t}_{i}$, where $\mathrm{in}_{i}$ stands for the canonical inclusion of $A_{i}$ into $\coprod_{i\in I} A_{i}$. 

It easily follows that $u^{t}$ is a homomorphism and satisfies the above equality.
\end{proof}

We next define the assignment $\mathfrak{L}_{\Sigma}$ from $\int^{\mathsf{NPIAlg}(\Sigma)_{\mathsf{s}}}\mathrm{Lsys}_{\Sigma}$ to 
$\int_{\mathsf{NPIAlg}(\Sigma)_{\mathsf{s}}}(\mathsf{Alg}(\Sigma)\,{\downarrow_{\mathsf{s}}}\, \bigcdot)$.

\begin{definition}\label{DLaSumFun} 
We let $\mathfrak{L}_{\Sigma}$ stand for the assignment from $\int^{\mathsf{NPIAlg}(\Sigma)_{\mathsf{s}}}\mathrm{Lsys}_{\Sigma}$ to 
$\int_{\mathsf{NPIAlg}(\Sigma)_{\mathsf{s}}}(\mathsf{Alg}(\Sigma)\,{\downarrow_{\mathsf{s}}}\, \bigcdot)$ defined as follows:
\begin{enumerate}
\item for every Lallement system of $\Sigma$-algebras $\mathbcal{A}=(\mathbf{I},\mathcal{A})$, $\mathfrak{L}_{\Sigma}(\mathbcal{A})$ is the pair $(\mathbf{I},(\boldsymbol{\mathfrak{L}}_{\Sigma}(\mathbcal{A}), \mathrm{p}_{\mathbcal{A}}))$, where $\boldsymbol{\mathfrak{L}}_{\Sigma}(\mathbcal{A})$ is the Lallement sum of $\mathbcal{A}$ as set forth in Definition~\ref{DLaSum}, and $\mathrm{p}_{\mathbcal{A}}\colon \boldsymbol{\mathfrak{L}}_{\Sigma}(\mathbcal{A})\mor \mathbf{I}$ the surjective homomorphism set forth in Proposition~\ref{PLaSum}, and
\item for every Lallement system of $\Sigma$-algebras $\mathbcal{A}=(\mathbf{I},\mathcal{A})$ and $\mathbcal{B}=(\mathbf{P},\mathcal{B})$ and every morphism $(t,u)\colon \mathbcal{A}\mor\mathbcal{B}$, $\mathfrak{L}_{\Sigma}(t,u)$ is the pair $(t,u^{t})$ where $u^{t}\colon \boldsymbol{\mathfrak{L}}_{\Sigma}(\mathbcal{A})\mor \boldsymbol{\mathfrak{L}}_{\Sigma}(\mathbcal{B})$ is the 
homomorphism set forth in Proposition~\ref{PLaSumMor}.
\end{enumerate}
\end{definition}

\begin{proposition}\label{PLaSumFun} 
The assignment $\mathfrak{L}_{\Sigma}$ is a covariant functor from $\int^{\mathsf{NPIAlg}(\Sigma)_{\mathsf{s}}}\mathrm{Lsys}_{\Sigma}$ to 
$\int_{\mathsf{NPIAlg}(\Sigma)_{\mathsf{s}}}(\mathsf{Alg}(\Sigma)\,{\downarrow_{\mathsf{s}}}\, \bigcdot)$. We will refer to $\mathfrak{L}_{\Sigma}$ as the \emph{functor of Lallement relative to 
$\Sigma$}.
\end{proposition}
\begin{proof}
That $\mathfrak{L}_{\Sigma}$ maps objects and morphisms of $\int^{\mathsf{NPIAlg}(\Sigma)_{\mathsf{s}}}\mathrm{Lsys}_{\Sigma}$ to objects and morphisms of
$\int_{\mathsf{NPIAlg}(\Sigma)_{\mathsf{s}}}(\mathsf{Alg}(\Sigma)\,{\downarrow_{\mathsf{s}}}\, \cdot)$ follows from Definition~\ref{DLaSum} and Propositions~\ref{PLaSum} and~\ref{PLaSumMor}, respectively. 

It is straightforward to prove that $\mathfrak{L}_{\Sigma}$ preserves identities and compositions.
\end{proof}

\section{Semi-inductive Lallement systems of $\Sigma$-algebras relative to an npi $\Sigma$-algebra and proof that for them the Lallement functor is a weak right multiadjoint}

In this section, after introducing the notion of semi-inductive Lallement systems of $\Sigma$-algebras relative to a naturally preordered idempotent $\Sigma$-algebra, we prove that the restriction of the functor $\mathfrak{L}_{\Sigma}$ to these Lallement systems is a weak right multiadjoint.

We begin by defining the notion of congruence on a $\Sigma$-algebra that preserves a subalgebra of it.

\begin{definition}\label{DEnv}
Let $\mathbf{A}$ be a subalgebra of $\mathbf{M}$ and $\Phi$ a congruence on $\mathbf{M}$. We will say that $\Phi$ \emph{preserves} $\mathbf{A}$ if $ \mathrm{pr}_{\Phi}\circ \mathrm{in}_{\mathbf{A},\mathbf{M}}$ is an injective homomorphism from $\mathbf{A}$ to $\mathbf{M}/\Phi$. 
\end{definition}

We next define, following~\cite{La67} and \cite{RS02}, the notion of semi-inductive Lallement system of 
$\Sigma$-algebras relative to an npi $\Sigma$-algebra.

\begin{definition}\label{DGCL}
Let $\mathcal{A}=((\mathbf{A}_{i}, \mathbf{M}_{i})_{i\in I}, (f_{i,j})_{(i,j)\in\leq_{\mathbf{I}}})$ be a Lallement system of $\Sigma$-algebras relative to an npi $\Sigma$-algebra $\mathbf{I}$. We will say that $\mathcal{A}$ is a \emph{semi-inductive Lallement system of $\Sigma$-algebras relative to $\mathbf{I}$}, or an \emph{$\mathbf{I}$-semi-inductive Lallement system of $\Sigma$-algebras}, if, 
\begin{enumerate} 
\item[(c)] for every $i\in I$ and every $a\in A_{i}$, $f_{i,i}(a) = a$, i.e., if $f_{i,i} = \mathrm{in}_{\mathbf{A}_{i},\mathbf{M}_{i}}$, the canonical embedding of $\mathbf{A}_{i}$ into $\mathbf{M}_{i}$.
\end{enumerate}
\end{definition}

\begin{remark}
In Definition~\ref{DGCL} we use the term ``semi-inductive'' because for $\mathcal{A}=((\mathbf{A}_{i}, \mathbf{M}_{i})_{i\in I}, (f_{i,j})_{(i,j)\in\leq_{\mathbf{I}}})$, a Lallement system of $\Sigma$-algebras relative to an npi $\Sigma$-algebra $\mathbf{I}$, and $i$, $j$, $k\in I$, with $(i,j)\in\leq_{\mathbf{I}}$ and $(j,k)\in\leq_{\mathbf{I}}$, we do not require that $f_{i,j}[A_{i}]\subseteq A_{j}$ and that, in addition, $f_{i,k} = f_{j,k}\circ f_{i,j}$; otherwise we would have a fully functorial definition of the concept in question, which is not the spirit either of~\cite{La67} or of~\cite{RS02}. Moreover, the condition $(c)$ in Definition~\ref{DGCL} is the condition $3)$ stated by Lallement in~\cite[Th\'{e}or\`{e}me 2.19, pp.~83--84]{La67}, as generalized by Romanowska and Smith in~\cite{RS91} to the field of universal algebra.
\end{remark}

\begin{assumption}\label{AGCL}
From now on, we will only consider semi-inductive Lallement system of $\Sigma$-algebras relative to an npi $\Sigma$-algebra and sums of them, and, occasionally, to shorten the writing, we will speak of ``Lallement systems'' instead of ``semi-inductive Lallement systems''.
Moreover, to simplify the notation, we continue to write (1) $\mathbf{I}$-$\mathsf{LAlg}(\Sigma)$ for the category of $\mathbf{I}$-semi-inductive Lallement systems; (2) $\mathrm{Lsys}_{\Sigma}$ for the  contravariant functor from $\mathsf{NPIAlg}(\Sigma)_{\mathsf{s}}$ to $\mathsf{Cat}$ that sends an npi $\Sigma$-algebra $\mathbf{I}$ to the category of $\mathbf{I}$-semi-inductive Lallement system of $\Sigma$-algebras; and (3) $\mathfrak{L}_{\Sigma}$ for the restriction of the Lallement functor 
$\mathfrak{L}_{\Sigma}$, set forth in Definition~\ref{PLaSumFun}, to the subcategory of $\int^{\mathsf{NPIAlg}(\Sigma)_{\mathsf{s}}}\mathrm{Lsys}_{\Sigma}$ determined by the semi-inductive Lallement systems labeled by npi $\Sigma$-algebras (which is also the domain of the split fibration  canonically associated to the contravariant functor mentioned in (2)), which shall also be referred to as the Lallement functor.
\end{assumption}

One of the most important features of semi-inductive Lallement systems is that their Lallement sums contain isomorphic copies of the subalgebras of their algebraic pairs.

\begin{proposition}\label{PGCL} 
Let $\mathcal{A}=((\mathbf{A}_{i}, \mathbf{M}_{i})_{i\in I}, (f_{i,j})_{(i,j)\in\leq_{\mathbf{I}}})$ be a semi-inductive Lallement system of $\Sigma$-algebras relative to an npi $\Sigma$-algebra $\mathbf{I}$, $\boldsymbol{\mathfrak{L}}_{\Sigma}(\mathbcal{A})$ its Lallement sum, set forth in Definition~\ref{DLaSum}, which, we recall, has as underlying set
$\coprod_{i\in I} A_{i}$.
Then, for every $i\in I$, $A_{i}\times\{i\}$ is a closed subset of $\boldsymbol{\mathfrak{L}}_{\Sigma}(\mathbcal{A})$ and we let $\mathbf{A}_{i}\times\{i\}$ stand for the $\Sigma$-subalgebra of 
$\boldsymbol{\mathfrak{L}}_{\Sigma}(\mathbcal{A})$ determined by $A_{i}\times\{i\}$. Moreover, the canonical projection onto the first component, i.e.,
\[
\pi^{i}_{0}\left\lbrace
\begin{array}{ccl}
\mathbf{A}_{i}\times\{i\}&\mor &\mathbf{A}_{i}\\
(a,i)&\longmapsto&a
\end{array}
\right.
\]
is an isomorphism.
\end{proposition}

Our next aim is to prove that the functor $\mathfrak{L}_{\Sigma}$ is a weak right multiadjoint. This will follow from the following series of propositions.

\begin{proposition}\label{PMult1} 
Let $(\mathbf{I},(\mathbf{A},f))$ be an object in $\int_{\mathsf{NPIAlg}(\Sigma)_{\mathsf{s}}}(\mathsf{Alg}(\Sigma)\,{\downarrow_{\mathsf{s}}}\, \bigcdot)$ and, for every $i\in I$, let $A_{i}$ be the fiber of $i$ along $f$, i.e., $A_{i}=f^{-1}[\{i\}]$. Then $A_{i}$ is a closed subset of the $\Sigma$-algebra $\mathbf{A}$ and we will let $\mathbf{A}_{i}$ stand for the subalgebra of $\mathbf{A}$ determined by $A_{i}$.
\end{proposition}

%
%

\begin{proposition}\label{PMult2} 
Let $(\mathbf{I},(\mathbf{A},f))$ be an object in $\int_{\mathsf{NPIAlg}(\Sigma)_{\mathsf{s}}}(\mathsf{Alg}(\Sigma)\,{\downarrow_{\mathsf{s}}}\, \bigcdot)$, $(\mathbf{A}_{i})_{i\in I}$ the family of $\Sigma$-algebras associated to the fibers of 
$\mathbf{A}$ and, for every $j\in I$, let $M_{j}$ stand for $\coprod_{i\in\Downarrow j} A_{i}$, where $\Downarrow\! j = \{i\in I\mid (i,j)\in \leq_{\mathbf{I}}\}$. Then there exists a structure of 
$\Sigma$-algebra $F^{\mathbf{M}_{j}}$ on $M_{j}$, we will denote by $\mathbf{M}_{j}=(M_{j}, F^{\mathbf{M}_{j}})$ the corresponding $\Sigma$-algebra, such that, for every $i,j\in I$, if $(i,j)\in\leq_{\mathbf{I}}$, then the canonical embedding $\mathrm{in}_{i,j}$ of $A_{i}$ into $M_{j}$ determines an injective homomorphism, also denoted by $\mathrm{in}_{i,j}$, of $\mathbf{A}_{i}$ into  
$\mathbf{M}_{j}$.
Thus, for every $i,j\in I$, if $(i,j)\in\leq_{\mathbf{I}}$, then $\mathbf{A}_{i}$ is isomorphic to the subalgebra $\mathrm{in}_{i,j}[\mathbf{A}_{i}]$ of $\mathbf{M}_{j}$. In particular, for every $i\in I$, $\mathbf{A}_{i}$ is isomorphic to the subalgebra $\mathrm{in}_{i,i}[\mathbf{A}_{i}]$ of $\mathbf{M}_{i}$.
\end{proposition}

\begin{proof}
Let $n$ be an element of $\mathbb{N}-\{0\}$, $\sigma\in \Sigma_{n}$ and $(a_{i_{k}}, i_{k})_{k\in n}$ a family of elements in $M_{j}$. Then, for every $k\in n$, $i_{k}\in \Downarrow\! j$ and $a_{i_{k}}\in A_{i_{k}}$. We define the structural operation associated to $\sigma$ acting on $(a_{i_{k}}, i_{k})_{k\in n}$ as 
\[
F^{\mathbf{M}_{j}}_{\sigma}((a_{i_{k}}, i_{k})_{k\in n})=
(
F^{\mathbf{A}}_{\sigma}((a_{i_{k}})_{k\in n}), F^{\mathbf{I}}_{\sigma}(({i_{k}})_{k\in n})
).
\]
In order to prove that this operation is well-defined we need to prove that (1) the element $i=F^{\mathbf{I}}_{\sigma}(({i_{k}})_{k\in n})$ is in $\Downarrow\! j$ and that (2) $F^{\mathbf{A}}_{\sigma}((a_{i_{k}})_{k\in n})\in A_{i}$.

Let us note that, for every $k\in n$, $i_{k}\in \Downarrow\! j$, i.e., $(i_{k},j)\in \leq_{\mathbf{I}}$.
Since $\mathbf{I}$ is naturally preordered, $F^{\mathbf{I}}_{\sigma}((i_{k})_{k\in n})\leq_{\mathbf{I}}F^{\mathbf{I}}_{\sigma}((j)_{k\in n})$. Now, since $\mathbf{I}$ is idempotent, we conclude that $F^{\mathbf{I}}_{\sigma}((j)_{k\in n})=j$. Hence, $i=F^{\mathbf{I}}_{\sigma}(({i_{k}})_{k\in n})\in \Downarrow\! j$. This proves (1).

Regarding (2), let us note that the following chain of equalities holds
\begin{align*}
f(F^{\mathbf{A}}_{\sigma}((a_{i_{k}})_{k\in n}))&=
F^{\mathbf{I}}_{\sigma}((f(a_{i_{k}}))_{k\in n})
\tag{$f$ is a homomorphism}
\\&=
F^{\mathbf{I}}_{\sigma}((i_{k})_{k\in n}).
\tag{$a_{i_{k}}\in A_{i_{k}}$}
\end{align*}
Now, taking into account that $i$ stands for $F^{\mathbf{I}}_{\sigma}((i_{k})_{k\in n})$, we infer that $F^{\mathbf{A}}_{\sigma}((a_{i_{k}})_{k\in n})\in f^{-1}[\{i\}]=A_{i}$. This proves (2).

All in all, we conclude that $\mathbf{M}_{j}=(M_{j}, F^{\mathbf{M}_{j}})$ is a well-defined $\Sigma$-algebra.
\end{proof}

\begin{definition}\label{DMult3} 
Let $(\mathbf{I},(\mathbf{A},f))$ be an object in $\int_{\mathsf{NPIAlg}(\Sigma)_{\mathsf{s}}}(\mathsf{Alg}(\Sigma)\,{\downarrow_{\mathsf{s}}}\, \bigcdot)$, $(\mathbf{A}_{i})_{i\in I}$ the family of $\Sigma$-algebras associated to the fibers of $\mathbf{A}$ set forth in Proposition~\ref{PMult1}, and $(\mathbf{M}_{i})_{i\in I}$ the family of $\Sigma$-algebras set forth in Proposition~\ref{PMult2}. Then we define the set of \emph{preserving congruences relative to $(\mathbf{I},(\mathbf{A},f))$} as 
\[
\textstyle
\mathrm{PCgr}(\mathbf{I},(\mathbf{A},f))=
\prod_{i\in I}
\mathrm{PCgr}(\mathbf{M}_{i}),
\]
where, for every $i\in I$, $\mathrm{PCgr}(\mathbf{M}_{i}) = \{\Phi\in \mathrm{Cgr}(\mathbf{M}_{i})\mid \Phi\mbox{ preserves }\mathrm{in}_{i,i}[\mathbf{A}_{i}]\}$.
\end{definition}

\begin{proposition}\label{PMult4} 
Let $(\mathbf{I},(\mathbf{A},f))$ be an object in $\int_{\mathsf{NPIAlg}(\Sigma)_{\mathsf{s}}}(\mathsf{Alg}(\Sigma)\,{\downarrow_{\mathsf{s}}}\, \bigcdot)$, $(\mathbf{A}_{i})_{i\in I}$ the family of $\Sigma$-algebras associated to the fibers of $\mathbf{A}$ set forth in Proposition~\ref{PMult1}, and $(\mathbf{M}_{i})_{i\in I}$ the family of $\Sigma$-algebras set forth in Proposition~\ref{PMult2}. Then, for every $\overline{\Theta} = (\Theta_{i})_{i\in I}\in\mathrm{PCgr}(\mathbf{I},(\mathbf{A},f))$ and every $i,j\in I$, if $(i,j)\in\leq_{\mathbf{I}}$, then the mapping $\mathrm{pr}_{\Theta_{j}}\circ \mathrm{in}_{i,j}\colon \mathbf{A}_{i}\mor \mathbf{M}_{j}/{\Theta_{j}}$ is a homomorphism.
Moreover, for every $i\in I$, $\mathbf{A}_{i}$ is isomorphic to the subalgebra $(\mathrm{pr}_{\Theta_{i}}\circ \mathrm{in}_{i,i})[\mathbf{A}_{i}]$ of $\mathbf{M}_{i}/{\Theta_{i}}$.
\end{proposition}

%
%

\begin{proposition}\label{PMult5} 
Let $(\mathbf{I},(\mathbf{A},f))$ be an object in $\int_{\mathsf{NPIAlg}(\Sigma)_{\mathsf{s}}}(\mathsf{Alg}(\Sigma)\,{\downarrow_{\mathsf{s}}}\, \bigcdot)$, $(\mathbf{A}_{i})_{i\in I}$ the family of $\Sigma$-algebras associated to the fibers of $\mathbf{A}$ set forth in Proposition~\ref{PMult1}, $(\mathbf{M}_{i})_{i\in I}$ the family of $\Sigma$-algebras set forth in Proposition~\ref{PMult2}, $\overline{\Theta}=(\Theta_{i})_{i\in I}\in\mathrm{PCgr}(\mathbf{I},(\mathbf{A},f))$ and, for every $i,j\in I$ such that $(i,j)\in\leq_{\mathbf{I}}$, let $g^{\overline{\Theta}}_{i,j}$ stand for the mapping from $\mathrm{p}_{\Theta_{i}}[\mathbf{A}_{i}]$ to $\mathbf{M}_{j}/\Theta_{j}$ defined as follows:
\[
g^{\overline{\Theta}}_{i,j}=\mathrm{pr}_{\Theta_{j}}\circ \mathrm{in}_{i,j}\circ 
\widehat{\mathrm{p}}_{\Theta_{i}}^{-1},
\]
where, for every $i\in I$, $\mathrm{p}_{\Theta_{i}}$ stands for $\mathrm{pr}_{\Theta_{i}}\circ \mathrm{in}_{i,i}$ and $\widehat{\mathrm{p}}_{\Theta_{i}}$ for the corestriction of $\mathrm{p}_{\Theta_{i}}$ to $\mathrm{p}_{\Theta_{i}}[\mathbf{A}_{i}]$. Then the ordered pair
\[
\mathcal{A}_{\overline{\Theta}}=
((\mathrm{p}_{\Theta_{i}}[\mathbf{A}_{i}], \mathbf{M}_{i}/{\Theta_{i}})_{i\in I}, (g^{\overline{\Theta}}_{i,j})_{(i,j)\in\leq_{\mathbf{I}}})
\]
is a semi-inductive Lallement system of $\Sigma$-algebras relative to $\mathbf{I}$.
\end{proposition}

\begin{proposition}\label{PMult6} 
Let $(\mathbf{I},(\mathbf{A},f))$ be an object in $\int_{\mathsf{NPIAlg}(\Sigma)_{\mathsf{s}}}(\mathsf{Alg}(\Sigma)\,{\downarrow_{\mathsf{s}}}\, \bigcdot)$, $\overline{\Theta}=(\Theta_{i})_{i\in I}\in\mathrm{PCgr}(\mathbf{I},(\mathbf{A},f))$ and $\mathbcal{A}_{\overline{\Theta}}=(\mathbf{I}, \mathcal{A}_{\overline{\Theta}})$ the semi-inductive Lallement system of $\Sigma$-algebras set forth in Proposition~\ref{PMult5}. Then there exists a morphism $(t_{\overline{\Theta}},h_{\overline{\Theta}})$ in $\int_{\mathsf{NPIAlg}(\Sigma)_{\mathsf{s}}}(\mathsf{Alg}(\Sigma)\,{\downarrow_{\mathsf{s}}}\, \bigcdot)$ from $(\mathbf{I},(\mathbf{A},f))$ to $\mathfrak{L}_{\Sigma}(\mathbcal{A}_{\overline{\Theta}})$. 
\end{proposition}

\begin{proof}
By Definition~\ref{DLaSumFun}, we have that $\mathfrak{L}_{\Sigma}(\mathbcal{A}_{\overline{\Theta}})=(\mathbf{I},(\boldsymbol{\mathfrak{L}}_{\Sigma}(\mathbcal{A}_{\overline{\Theta}}), \mathrm{p}_{\mathbcal{A}_{\overline{\Theta}}}))$. We define $t_{\overline{\Theta}} = \mathrm{id}_{\mathbf{I}}$, which, in particular, is a surjective homomorphism from $\mathbf{I}$ to itself.

On the other hand, by Proposition~\ref{PMult5}, we have that 
\[
\textstyle
\mathcal{A}_{\overline{\Theta}}=
((\mathrm{p}_{\Theta_{i}}[\mathbf{A}_{i}], \mathbf{M}_{i}/{\Theta_{i}})_{i\in I}, (g^{\overline{\Theta}}_{i,j})_{(i,j)\in\leq_{\mathbf{I}}})
\]
By Definition~\ref{DLaSum}, we have that the underlying set of $\boldsymbol{\mathfrak{L}}_{\Sigma}(\mathbcal{A}_{\overline{\Theta}})$ is 
\[
\textstyle
\coprod_{i\in I}\mathrm{p}_{\Theta_{i}}[A_{i}]=\bigcup_{i\in I}(\mathrm{p}_{\Theta_{i}}[A_{i}]\times \{i\}).
\]
We define the mapping $h_{\overline{\Theta}}$ as follows
\[
h_{\overline{\Theta}}
\left\lbrace
\begin{array}{ccl}
A&\mor&\coprod_{i\in I}\mathrm{p}_{\Theta_{i}}[A_{i}]\\
a&\longmapsto& (\mathrm{p}_{\Theta_{f(a)}}(a),f(a))
\end{array}
\right.
\]

One can easily prove that the mapping $h_{\overline{\Theta}}$ determines a homomorphism from $\mathbf{A}$ to  $\boldsymbol{\mathfrak{L}}_{\Sigma}(\mathbcal{A}_{\overline{\Theta}})$ such that $t_{\overline{\Theta}}\circ f=\mathrm{p}_{\mathbcal{A}_{\overline{\Theta}}}\circ h_{\overline{\Theta}}$.
\end{proof}

The following result is the centerpiece of this work and the one from which its title is derived. Here, we establish that the functor $\mathfrak{L}_{\Sigma}$ is a weak right multiadjoint.

\begin{theorem}\label{PMult7} 
Let $(\mathbf{I},(\mathbf{A},f))$ be an object in 
$\int_{\mathsf{NPIAlg}(\Sigma)_{\mathsf{s}}}(\mathsf{Alg}(\Sigma)\,{\downarrow_{\mathsf{s}}}\, \bigcdot)$, $\mathbcal{B}=(\mathbf{P},\mathcal{B})$ a semi-inductive Lallement system of $\Sigma$-algebras in $\int^{\mathsf{NPIAlg}(\Sigma)_{\mathsf{s}}}\mathrm{Lsys}_{\Sigma}$ and $(t,h)$ a morphism in 
$\int_{\mathsf{NPIAlg}(\Sigma)_{\mathsf{s}}}(\mathsf{Alg}(\Sigma)\,{\downarrow_{\mathsf{s}}}\, \bigcdot)$ from $(\mathbf{I},(\mathbf{A},f))$ to $\mathfrak{L}_{\Sigma}(\mathbcal{B})$. Then there exists a non-empty subset $\mathrm{PCgr}(\mathbcal{B})\subseteq \mathrm{PCgr}(\mathbf{I},(\mathbf{A},f))$ such that, for every $\overline{\Theta}=(\Theta_{i})_{i\in I}\in \mathrm{PCgr}(\mathbcal{B})$, there exists a unique morphism $(s_{\overline{\Theta}}, v_{\overline{\Theta}})$ in $\int^{\mathsf{NPIAlg}(\Sigma)_{\mathsf{s}}}\mathrm{Lsys}_{\Sigma}$ from $\mathbcal{A}_{\overline{\Theta}}$ to $\mathbcal{B}$ such that 
\[
(t,h)=\mathfrak{L}_{\Sigma}(s_{\overline{\Theta}}, v_{\overline{\Theta}})\circ (t_{\overline{\Theta}}, h_{\overline{\Theta}}),
\]
where $\mathbcal{A}_{\overline{\Theta}}=(\mathbf{I}, \mathcal{A}_{\overline{\Theta}})$ is the semi-inductive Lallement system of $\Sigma$-algebras set forth in Proposition~\ref{PMult5} and 
$
(t_{\overline{\Theta}}, h_{\overline{\Theta}})\colon
(\mathbf{I},(\mathbf{A},f))
\mor
\mathfrak{L}_{\Sigma}(\mathbcal{A}_{\overline{\Theta}})
$
the morphism in $\int_{\mathsf{NPIAlg}(\Sigma)_{\mathsf{s}}}(\mathsf{Alg}(\Sigma)\,{\downarrow_{\mathsf{s}}}\, \bigcdot)$ set forth in Proposition~\ref{PMult6}. 
\end{theorem}

\begin{proof}
\textsf{Existence.}  By Definition~\ref{DLaSumFun}, we have that $\mathfrak{L}_{\Sigma}(\mathbcal{B})=(\mathbf{P},(\boldsymbol{\mathfrak{L}}_{\Sigma}(\mathbcal{B}), \mathrm{p}_{\mathbcal{B}}))$, where $\boldsymbol{\mathfrak{L}}_{\Sigma}(\mathbcal{B})$ is the Lallement sum of $\mathbcal{B}$, set forth in Definition~\ref{DLaSum}, and $\mathrm{p}_{\mathbcal{B}}$ 
the surjective homomorphism from $\boldsymbol{\mathfrak{L}}_{\Sigma}(\mathbcal{B})$ to $\mathbf{P}$
set forth in Proposition~\ref{PLaSum}. By Definition~\ref{DAlgdGC}, the morphism $(t,h)\colon (\mathbf{I},(\mathbf{A},f))\mor \mathfrak{L}_{\Sigma}(\mathbcal{B})$ in $\int_{\mathsf{NPIAlg}(\Sigma)_{\mathsf{s}}}(\mathsf{Alg}(\Sigma)\,{\downarrow_{\mathsf{s}}}\, \bigcdot)$ is given by a surjective homomorphism $t\colon \mathbf{I}\mor \mathbf{P}$ together with a homomorphism $h\colon \mathbf{A}\mor \boldsymbol{\mathfrak{L}}_{\Sigma}(\mathbcal{B})$ such that 
\begin{equation}
t\circ f=\mathrm{p}_{\mathbcal{B}}\circ h.
\tag{E1}\label{E1}
\end{equation}

Let us assume that $\mathcal{B}=((\mathbf{B}_{p}, \mathbf{N}_{p})_{p\in P}, (g_{p,q})_{(p,q)\in \leq_{\mathbf{P}}})$. By Definition~\ref{DLaSum}, the underlying set of $\boldsymbol{\mathfrak{L}}_{\Sigma}(\mathbcal{B})$ is given by
\[
\textstyle
\mathfrak{L}_{\Sigma}(\mathbcal{B})=\coprod_{p\in P}B_{p}=\bigcup_{p\in P}(B_{p}\times \{p\}).
\]

By Proposition~\ref{PMult1}, let $(\mathbf{A}_{i})_{i\in I}$ be the family of $\Sigma$-algebras associated to the fibers of $\mathbf{A}$.

\begin{claim}\label{CMul71} 
Let $a$ be an element of $A$. Then $h(a)\in \mathrm{B}_{t(f(a))}\times \{t(f(a))\}$. In particular, for every $j\in I$, we have that $h[A_{j}]\subseteq B_{t(j)}\times \{t(j)\}$. 
\end{claim}
For every element $a$ of $A$ we have that
\allowdisplaybreaks
\begin{align*}
\mathrm{p}_{\mathbcal{B}}(h(a))&=
t(f(a)).
\tag{by Eq.~\ref{E1}}
\end{align*}
From the last equality, we conclude that $h(a)\in \mathrm{p}_{\mathbcal{B}}^{-1}[\{t(f(a))\}]$. Let us recall that, by Proposition~\ref{PLaSum}, $\mathrm{p}_{\mathbcal{B}}$ is the homomorphism given by
\[
\mathrm{p}_{\mathbcal{B}}\colon
\left\lbrace
\begin{array}{rcl}
\boldsymbol{\mathfrak{L}}_{\Sigma}(\mathbcal{B})&\mor&\mathbf{P}\\
(b,p)&\longmapsto&p
\end{array}
\right.
\]
Thus, $\mathrm{p}_{\mathbcal{B}}^{-1}[\{t(f(a))\}]=B_{t(f(a))}\times \{t(f(a))\}$.   In particular, for every $j\in I$ and every $a\in A$, since $a\in A_{j}$ if, and only if, $f(a)=j$, we have that $h[A_{j}]\subseteq B_{t(j)}\times \{t(j)\}$.  

This completes the proof of Claim~\ref{CMul71}.

On the basis of Claim~\ref{CMul71}, for every $j\in I$, we let $h_{j}$ stand for the birestriction of $h$ to $A_{j}$ and $B_{t(j)}\times \{t(j)\}$. Note that $h_{j}$ determines a homomorphism from $\mathbf{A}_{j}$ to $\mathbf{B}_{t(j)}\times \{t(j)\}$ (see Definition~\ref{DSubAlg}). Moreover, we denote by $\pi^{t(j)}_{0}$ the canonical projection from $B_{t(j)}\times \{t(j)\}$ to $B_{t(j)}$, which, by Proposition~\ref{PGCL}, is an isomorphism from $\mathbf{B}_{t(j)}\times \{t(j)\}$ to $\mathbf{B}_{t(j)}$. Then  
$\pi^{t(j)}_{0}\circ h_{j}$ is a homomorphism from  $\mathbf{A}_{j}$ to $\mathbf{B}_{t(j)}$.

Let $i$ and $j$ be elements of $I$ such that $(i,j)\in\leq_{\mathbf{I}}$. Then, by  Proposition~\ref{PSubOrd}, we have that $(t(i),t(j))\in\leq_{\mathbf{P}}$ and we let $r_{i,j}$ stand for the mapping from $A_{i}$ to $N_{t(j)}$ defined as follows:
\begin{align*}
\textstyle
r_{i,j} = g_{t(i),t(j)}\circ \pi^{t(i)}_{0}\circ h_{i}.
\tag{E2}\label{E2}
\end{align*}
Then, since $g_{t(i),t(j)}$, $\pi^{t(i)}_{0}$ and  $h_{i}$ are homomorphisms we can assert that $r_{i,j}$ is a homomorphism from $\mathbf{A}_{i}$ to $\mathbf{N}_{t(j)}$.

Since it will be used several times later, we consider, for $i\in I$, the particular case of the homomorphism $r_{i,i}$ from $\mathbf{A}_{i}$ to $\mathbf{N}_{t(i)}$.
\begin{claim}\label{CMul72} Let $i$ be an element of $I$. Then $r_{i,i}=\pi^{t(i)}_{0}\circ h_{i}$. Hence, for every $a\in A_{i}$, $r_{i,i}(a)=  \pi^{t(i)}_{0}(h(a))$.
\end{claim}
This follows from the definition of $r_{i,j}$ in Equation~\ref{E2} for the case in which $i=j$, Definition~\ref{DGCL} and the definition of $h_{i}$ in~Claim~\ref{CMul71}.

This completes the proof of Claim~\ref{CMul72}.

Let $j$ be an element of $I$. Then, by Proposition~\ref{PMult2}, we have that $M_{j}=\coprod_{i\in\Downarrow j} A_{i}$, where $\Downarrow\! j = \{i\in I\mid (i,j)\in \leq_{\mathbf{I}}\}$. Hence, by the universal property of the coproduct of sets, there exists a unique mapping $[(r_{i,j})_{i\in\Downarrow j}]\colon M_{j}\mor N_{t(j)}$, denoted by $r_{j}$ for short, such that, for every $i\in \Downarrow j$,
\begin{align*}
r_{j}\circ \mathrm{in}_{i,j}=r_{i,j},
\tag{E3}\label{E3}
\end{align*}
where $\mathrm{in}_{i,j}$ stands for the canonical inclusion of $A_{i}$ into $M_{j}$, see Proposition~\ref{PMult2}.

\begin{claim}\label{CMul73} 
Let $j$ be an element of $I$. Then $r_{j}$ is a homomorphism from $\mathbf{M}_{j}$ to $\mathbf{N}_{t(j)}$.
\end{claim}
Let $n$ be an element of $\mathbb{N}-\{0\}$, $\sigma\in \Sigma_{n}$ and $(a_{k},i_{k})_{k\in n}$ a family of elements in $M_{j}^{n}$. In what follows, to simplify notation we let $i$ stand for $F^{\mathbf{I}}_{\sigma}(({i_{k}})_{k\in n})$. Then the following chain of equalities holds
\begin{flushleft}
$r_{j}(
F^{\mathbf{M}_{j}}_{\sigma}((a_{k},i_{k})_{k\in n})
)$
\allowdisplaybreaks
\begin{align*}
&=
r_{j}(
F^{\mathbf{A}}_{\sigma}((a_{i_{k}})_{k\in n}), i
)
\tag{Def.~$\mathbf{M}_{j}$, Prop.~\ref{PMult2}}
\\&=
r_{j}(\mathrm{in}_{i,j}(
F^{\mathbf{A}}_{\sigma}((a_{i_{k}})_{k\in n})
))
\tag{Def.~$\mathrm{in}_{i,j}$, Prop.~\ref{PMult2}}
\\&=
r_{i,j}(
F^{\mathbf{A}}_{\sigma}((a_{i_{k}})_{k\in n})
)
\tag{By Eq.~\ref{E3}}
\\&=
g_{t(i),t(j)}( 
\pi^{t(i)}_{0}( 
h_{i}(
F^{\mathbf{A}}_{\sigma}((a_{i_{k}})_{k\in n})
)))
\tag{By Eq.~\ref{E2}}
\\&=
g_{t(i),t(j)}( 
\pi^{t(i)}_{0}( 
h(
F^{\mathbf{A}}_{\sigma}((a_{i_{k}})_{k\in n})
)))
\tag{Def.~$h_{i}$} 
\\&=
g_{t(i),t(j)}( 
\pi^{t(i)}_{0}( 
F^{\boldsymbol{\mathfrak{L}}_{\Sigma}(\mathbcal{B})}_{\sigma}(
(h(a_{i_{k}}))_{k\in n})
))
\tag{$h$ is $\Sigma$-hom.} 
\\&=
g_{t(i),t(j)}( 
\pi^{t(i)}_{0}( 
F^{\boldsymbol{\mathfrak{L}}_{\Sigma}(\mathbcal{B})}_{\sigma}(
(
\pi^{t(i_{k})}_{0}(h(a_{i_{k}})),
t(i_{k})
)_{k\in n})
))
\tag{Claim.~\ref{CMul71}} 
\\&=
g_{t(i),t(j)}( 
\pi^{t(i)}_{0}( 
F^{\mathbf{N}_{t(i)}}_{\sigma}(
(g_{t(i_{k}),t(i)}(
\pi^{t(i_{k})}_{0}(h(a_{i_{k}}))
))_{k\in n}
),
t(i)
))
\tag{Def.~\ref{DLaSum}}
\\&=
g_{t(i),t(j)}( 
F^{\mathbf{N}_{t(i)}}_{\sigma}(
(g_{t(i_{k}),t(i)}(
\pi^{t(i_{k})}_{0}(h(a_{i_{k}}))
))_{k\in n}
)
)
\tag{Def.~$\pi^{t(i)}_{0}$}
\\&=
F^{\mathbf{N}_{t(j)}}_{\sigma}(
(
g_{t(i_{k}),t(j)}(
\pi^{t(i_{k})}_{0}(
h(
a_{i_{k}}
)
)
)
)_{k\in n}
)
\tag{Def.~\ref{DLsys}}
\\&=
F^{\mathbf{N}_{t(j)}}_{\sigma}(
(
g_{t(i_{k}),t(j)}(
\pi^{t(i_{k})}_{0}(
h_{i_{k}}(
a_{i_{k}}
)
)
)
)_{k\in n}
)
\tag{Def.~$h_{i_{k}}$} 
\\&=
F^{\mathbf{N}_{t(j)}}_{\sigma}(
(r_{i_{k},j}(a_{i_{k}}))_{k\in n}
)
\tag{By Eq.~\ref{E2}}
\\&=
F^{\mathbf{N}_{t(j)}}_{\sigma}(
(r_{j}(\mathrm{in}_{i_{k},j}(a_{i_{k}})))_{k\in n}
)
\tag{Def.~$\mathrm{in}_{i_{k},j}$, Prop.~\ref{PMult2}}
\\&=
F^{\mathbf{N}_{t(j)}}_{\sigma}(
(r_{j}(a_{i_{k}},i_{k}))_{k\in n}
).
\tag{By Eq.~\ref{E3}}
\end{align*}
\end{flushleft}

This completes the proof of Claim~\ref{CMul73}.

By Claim~\ref{CMul73} we have that, for every $i\in I$, $\mathrm{Ker}(r_{i})$ is a congruence on  
$\mathbf{M}_{i}$. Then we let $\mathrm{PCgr}(\mathbcal{B})$ stand for 
\[
\textstyle
\mathrm{PCgr}(\mathbcal{B})=
\prod_{i\in I}\{\Phi\in [\Delta_{\mathbf{M}_{i}}, \mathrm{Ker}(r_{i})]\mid \Phi\mbox{ preserves }\mathrm{in}_{i,i}[\mathbf{A}_{i}]\},
\]
where $[\Delta_{\mathbf{M}_{i}}, \mathrm{Ker}(r_{i})]$ is the closed interval of $\mathbf{Cgr}(\mathbf{M}_{i})$ determined by $\Delta_{\mathbf{M}_{i}}$ and $\mathrm{Ker}(r_{i})$. Let us note that $\mathrm{PCgr}(\mathbcal{B})$ is non-empty, since it contains $(\Delta_{\mathbf{M}_{i}})_{i \in I}$. Moreover, by Definition~\ref{DMult3}, $\mathrm{PCgr}(\mathbcal{B})$ is a subset of $\mathrm{PCgr}(\mathbf{I},(\mathbf{A},f))$. 

Let  $\overline{\Theta}=(\Theta_{i})_{i\in I}$ be an element in $\mathrm{PCgr}(\mathbcal{B})$. Then, by Proposition~\ref{PMult5}, we have the semi-inductive Lallement system of $\Sigma$-algebras $\mathbcal{A}_{\overline{\Theta}}=(\mathbf{I}, \mathcal{A}_{\overline{\Theta}})$, where
\[
\mathcal{A}_{\overline{\Theta}}=
((\mathrm{p}_{{\Theta}_{i}}[\mathbf{A}_{i}], \mathbf{M}_{i}/{{\Theta}_{i}})_{i\in I}, (g^{\overline{\Theta}}_{i,j})_{(i,j)\in\leq_{\mathbf{I}}}).
\]

We can now provide the morphism $(s_{\overline{\Theta}},v_{\overline{\Theta}})$ in $\int^{\mathsf{NPIAlg}(\Sigma)_{\mathsf{s}}}\mathrm{Lsys}_{\Sigma}$ from $\mathbcal{A}_{\overline{\Theta}}$ to $\mathbcal{B}=(\mathbf{P},\mathcal{B})$ such that $(t,h)=\mathfrak{L}_{\Sigma}(s_{\overline{\Theta}}, v_{\overline{\Theta}})\circ (t_{\overline{\Theta}}, h_{\overline{\Theta}})$. We define $s_{\overline{\Theta}}=t$. This is, by assumption, a surjective homomorphism from $\mathbf{I}$ to $\mathbf{P}$. On the other hand, by Definition~\ref{DLsysGC2}, $v_{\overline{\Theta}}=(v_{\overline{\Theta}, i})_{i\in I}$ must be a morphism of $\mathbf{I}$-Lallement systems of $\Sigma$-algebras from $\mathcal{A}_{\overline{\Theta}}$ to $\mathcal{B}_{t}$, where, by Definition~\ref{DLsysFun}, we must have that 
\[
\mathcal{B}_{t}=((\mathbf{B}_{t(i)}, \mathbf{N}_{t(i)})_{i\in I}, (g_{t(i),t(j)})_{(i,j)\in \leq_{\mathbf{I}}}).
\]

Let $i$ be an element of $I$. Then we define $v_{\overline{\Theta}, i}$ to be the unique mapping from $M_{i}/{\Theta_{i}}$ to $N_{t(i)}$ such that 
\begin{align*}
v_{\overline{\Theta}, i}\circ \mathrm{pr}_{\Theta_{i}}=r_{i}.
\label{E4}\tag{E4}
\end{align*}

This mapping exists and is well-defined because, for every $i\in I$, $\Theta_{i}\subseteq \mathrm{Ker}(r_{i})$.

\begin{claim}\label{CMul74} 
Let $j$ be an element of $I$. Then $v_{\overline{\Theta}, j}$ is a homomorphism from 
$\mathbf{M}_{j}/{\Theta_{j}}$ to $\mathbf{N}_{t(j)}$.
\end{claim}

This follows from the definition of $v_{\overline{\Theta}, j}$ and from the fact that $r_{i}$ and $\mathrm{pr}_{\Theta_{i}}$ are homomorphisms.

%
This completes the proof of Claim~\ref{CMul74}.

\begin{claim}\label{CMul75} 
Let $i$ be an element of $I$. Then $v_{\overline{\Theta}, i}[ \mathrm{p}_{\Theta_{i}}[A_{i}] ]\subseteq B_{t(i)}$.
\end{claim}
Let $i$ be an element of $I$ and $a\in A_{i}$. Then the following chain of equalities holds.
\allowdisplaybreaks
\begin{align*}
v_{\overline{\Theta}, i}(\mathrm{p}_{\Theta_{i}}(a))&=
v_{\overline{\Theta}, i}(\mathrm{pr}_{\Theta_{i}}(\mathrm{in}_{i,i}(a)))
\tag{Def.~$\mathrm{p}_{\Theta_{i}}$, Prop.~\ref{PMult5}}
\\&=
r_{i}(\mathrm{in}_{i,i}(a))
\tag{By Eq.~\ref{E4}}
\\&=
r_{i,i}(a)
\tag{By Eq.~\ref{E3}}
\\&=
\pi^{t(i)}_{0}(h(a)).
\tag{Claim~\ref{CMul72}}
\end{align*}

By Claim~\ref{CMul71}, $\pi^{t(i)}_{0}(h(a))\in B_{t(i)}$. Hence, $v_{\overline{\Theta}, i}[ \mathrm{p}_{\Theta_{i}}[A_{i}] ]\subseteq B_{t(i)}$.

This completes the proof of Claim~\ref{CMul75}.

\begin{claim}\label{CMul76} 
Let $i$ and $j$ be elements of $I$ such that $(i,j)\in\leq_{\mathbf{I}}$. Then $v_{\overline{\Theta},  j}\circ g^{\overline{\Theta}}_{i,j}=g_{t(i),t(j)}\circ v_{\overline{\Theta}, i}$.
\end{claim}

The following chain of equalities holds
\allowdisplaybreaks
\begin{align*}
v_{\overline{\Theta},  j}(g^{\overline{\Theta}}_{i,j}(
\mathrm{p}_{\Theta_{i}}(a)))
&=
v_{\overline{\Theta},  j}(g^{\overline{\Theta}}_{i,j}(
\mathrm{pr}_{\Theta_{i}}(\mathrm{in}_{i,i}(
a))
))
\tag{Def.~$\mathrm{p}_{\Theta_{i}}$, Prop.~\ref{PMult5}}
\\&=
v_{\overline{\Theta},  j}(g^{\overline{\Theta}}_{i,j}(
\mathrm{pr}_{\Theta_{i}}(a,i)
))
\tag{Def.~$\mathrm{in}_{i,i}$, Prop.~\ref{PMult2}}
\\&=
v_{\overline{\Theta},  j}(g^{\overline{\Theta}}_{i,j}(
[(a,i)]_{\Theta_{i}}
))
\tag{Def.~$\mathrm{pr}_{\Theta_{i}}$}
\\&=
v_{\overline{\Theta},  j}(
[(a,i)]_{\Theta_{j}}
)
\tag{Prop.~\ref{PMult5}}
\\&=
v_{\overline{\Theta},  j}(
\mathrm{pr}_{\Theta_{j}}(a,i)
)
\tag{Def.~$\mathrm{pr}_{\Theta_{j}}$}
\\&=
v_{\overline{\Theta},  j}(
\mathrm{pr}_{\Theta_{j}}(
\mathrm{in}_{i,j}(a)
)
)
\tag{Def.~$\mathrm{in}_{i,j}$, Prop.~\ref{PMult2}}
\\&=
r_{j}(
\mathrm{in}_{i,j}(a)
)
\tag{By Eq.~\ref{E4}}
\\&=
r_{i,j}(
a
)
\tag{By Eq.~\ref{E3}}
\\&=
g_{t(i),t(j)}(\pi^{t(i)}_{0}(h_{i}(a)))
\tag{By Eq.~\ref{E2}}
\\&=
g_{t(i),t(j)}(\pi^{t(i)}_{0}(h(a)))
\tag{Def.~$h_{i}$}
\\&=
g_{t(i),t(j)}(
v_{\overline{\Theta}, i}(
\mathrm{p}_{\Theta_{i}}(a)
)
).
\tag{Claim~\ref{CMul75}}
\end{align*}

This completes the proof of Claim~\ref{CMul76}.

\begin{claim}\label{CMul77} The family $v_{\overline{\Theta}}=(v_{\overline{\Theta}, i})_{i\in I}$   is a morphism of $\mathbf{I}$-Lallement systems of $\Sigma$-algebras $v_{\overline{\Theta}}\colon\mathcal{A}_{\overline{\Theta}}\mor  \mathcal{B}_{t}$.
\end{claim}
It follows from Definition~\ref{DLsys} and Claims~\ref{CMul74},~\ref{CMul75} and~\ref{CMul76}.

This completes the proof of Claim~\ref{CMul77}.

\begin{claim}\label{CMul78} 
Let $(v_{\overline{\Theta}})^{t}$ be the homomorphism introduced in Proposition~\ref{PLaSumMor}.
Then the equality $(v_{\overline{\Theta}})^{t}\circ  h_{\overline{\Theta}}=h$ holds.
\end{claim}

Let $a$ be an element of $A$. Then the following chain of equalities holds
\allowdisplaybreaks
\begin{align*}
(v_{\overline{\Theta}})^{t}(  h_{\overline{\Theta}}(a))&=
(v_{\overline{\Theta}})^{t}(\mathrm{p}_{\Theta_{f(a)}}(a),f(a))
\tag{Def.~$h_{\overline{\Theta}}$, Prop.~\ref{PMult6}}
\\&=
(v_{\overline{\Theta}})^{t}(\mathrm{in}_{f(a)}(
\mathrm{p}_{\Theta_{f(a)}}(a)
))
\tag{Def.~$\mathrm{in}_{f(a)}$}
\\&=
(v_{\overline{\Theta}})^{t}_{f(a)}(
\mathrm{p}_{\Theta_{f(a)}}(a)
)
\tag{Def.~$(v_{\overline{\Theta}})^{t}$, Prop.~\ref{PLaSumMor}}
\\&=
\mathrm{in}_{t(f(a))}(v_{\overline{\Theta}, t(a)}(
\mathrm{p}_{\Theta_{f(a)}}(a)
))
\tag{Def.~$(v_{\overline{\Theta}})^{t}_{f(a)}$, Prop.~\ref{PLaSumMor}}
\\&=
\mathrm{in}_{t(f(a))}(
\pi^{t(f(a))}_{0}(
h(a)
))
\tag{Claim~\ref{CMul75}}
\\&=
(\pi^{t(f(a))}_{0}(
h(a)
),
t(f(a)))
\tag{Def.~$\mathrm{in}_{t(f(a))}$}
\\&=
(\pi^{\mathrm{p}_{\mathbcal{B}}( h(a))}_{0}(
h(a)
),
\mathrm{p}_{\mathbcal{B}}( h(a))
)
\tag{By Eq.~\ref{E1}}
\\&=
h(a).
\tag{Def.~$\mathrm{p}_{\mathbcal{B}}$, Prop.~\ref{PLaSum}}
\end{align*}

This completes the proof of Claim~\ref{CMul78}

\begin{claim}\label{CMul79} The  equality $(t,h)=\mathfrak{L}_{\Sigma}(s_{\overline{\Theta}}, v_{\overline{\Theta}})\circ (t_{\overline{\Theta}}, h_{\overline{\Theta}})$ holds.
\end{claim}

The following chain of equalities holds
\begin{flushleft}
$\mathfrak{L}_{\Sigma}(s_{\overline{\Theta}}, v_{\overline{\Theta}})\circ (t_{\overline{\Theta}}, h_{\overline{\Theta}})$
\allowdisplaybreaks
\begin{align*}
&=
(s_{\overline{\Theta}}, (v_{\overline{\Theta}})^{s_{\overline{\Theta}}})\circ 
(t_{\overline{\Theta}}, h_{\overline{\Theta}})
\tag{Def.~$\mathfrak{L}_{\Sigma}$, Def.~\ref{DLaSumFun}}
\\&=
(t, (v_{\overline{\Theta}})^{t})\circ 
(t_{\overline{\Theta}}, h_{\overline{\Theta}})
\tag{Def.~$s_{\overline{\Theta}}$}
\\&=
(t, (v_{\overline{\Theta}})^{t})\circ 
(\mathrm{id}_{\mathbf{I}}, h_{\overline{\Theta}})
\tag{Def.~$t_{\overline{\Theta}}$, Prop.~\ref{PMult6}}
\\&=
(
t\circ \mathrm{id}_{\mathbf{I}}
, (v_{\overline{\Theta}})^{t}\circ  (\mathsf{Alg}(\Sigma)\,{\downarrow_{\mathsf{s}}}\, t)(h_{\overline{\Theta}}))
\tag{Composition, Def.~\ref{DGConsOp}}
\\&=
(
t
, (v_{\overline{\Theta}})^{t}\circ  (\mathsf{Alg}(\Sigma)\,{\downarrow_{\mathsf{s}}}\, t)(h_{\overline{\Theta}}))
\tag{Identity $\mathrm{id}_{\mathbf{I}}$}
\\&=
(
t
, (v_{\overline{\Theta}})^{t}\circ h_{\overline{\Theta}})
\tag{Def.~$(\mathsf{Alg}(\Sigma)\,{\downarrow_{\mathsf{s}}}\, t)$, Def.~\ref{DAlgda}}
\\&=
(t,h).
\tag{Claim~\ref{CMul79}}
\end{align*}
\end{flushleft}

This completes the proof of Claim~\ref{CMul79}.

This proves the existence of the desired morphism.

\textsf{Unicity.} Let $(s, v)\colon\mathbcal{A}_{\overline{\Theta}}\mor \mathbcal{B}$ be a morphism in $\int^{\mathsf{NPIAlg}(\Sigma)_{\mathsf{s}}}\mathrm{Lsys}_{\Sigma}$ such that 
\begin{align*}
(t,h)&=\mathfrak{L}_{\Sigma}(s, v)\circ (t_{\overline{\Theta}}, h_{\overline{\Theta}}).
\tag{U1}\label{U1}
\end{align*}

We recall that, by Definition~\ref{DLsysGC2}, the morphism $(s, v)$ is such that $s\colon \mathbf{I}\mor \mathbf{P}$ is a surjective homomorphism and $v=(v_{i})_{i\in I}$ a morphism of $\mathbf{I}$-Lallement systems from $\mathbcal{A}_{\overline{\Theta}}$ to $\mathbcal{B}_{s}$. Thus, by  Definition~\ref{DLsys}, for every $i\in I$, $v_{i}\colon \mathbf{M}_{i}/{\Theta_{i}}\mor \mathbf{N}_{s(i)}$ is a homomorphism such that $
v_{i}[\mathrm{p}_{\Theta_{i}}[\mathrm{A}_{i}]]\subseteq B_{s(i)}$. Moreover, for every $(i,j)\in\leq_{\mathbf{I}}$, the following equality holds
\begin{align*}
v_{j}\circ g^{\overline{\Theta}}_{i,j}&=g_{s(i),s(j)}\circ v_{i}.
\tag{U2}\label{U2}
\end{align*}
In addition, the following chain of equalities holds
\allowdisplaybreaks
\begin{align*}
(t,h)&=\mathfrak{L}_{\Sigma}(s, v)\circ (t_{\overline{\Theta}}, h_{\overline{\Theta}})
\tag{by Eq.~\ref{U1}}
\\&=
(s,v^{s})\circ (t_{\overline{\Theta}}, h_{\overline{\Theta}})
\tag{Def.~$\mathfrak{L}_{\Sigma}$, Def.~\ref{DLaSumFun}}
\\&=
(s,v^{s})\circ (\mathrm{id}_{\mathbf{I}}, h_{\overline{\Theta}})
\tag{Def.~$t_{\overline{\Theta}}$, Prop.~\ref{PMult6}}
\\&=
(s\circ \mathrm{id}_{\mathbf{I}}, v^{s}\circ  (\mathsf{Alg}(\Sigma)\,{\downarrow_{\mathsf{s}}}\, s)(h_{\overline{\Theta}}))
\tag{Composition, Def.~\ref{DGConsOp}}
\\&=
(s, v^{s}\circ  (\mathsf{Alg}(\Sigma)\,{\downarrow_{\mathsf{s}}}\, s)(h_{\overline{\Theta}}))
\tag{Identity ~$\mathrm{id}_{\mathbf{I}}$}
\\&=
(s, v^{s}\circ  h_{\overline{\Theta}}).
\tag{Def.~$(\mathsf{Alg}(\Sigma)\,{\downarrow_{\mathsf{s}}}\, s)$, Def.~\ref{DAlgda}}
\end{align*}

From the above equality we conclude that $s=t$ and that 
\begin{align*}
v^{s}\circ  h_{\overline{\Theta}}=v^{t}\circ  h_{\overline{\Theta}}=h.
\tag{U3}\label{U3}
\end{align*}

Taking into account Proposition~\ref{PLaSumMor}, we set, for every $i\in I$, $v^{t}_{i}=\mathrm{in}_{t(i)}\circ v_{i}$ and, finally, we set $v^{t}=[(v^{t}_{i})_{i\in I}]$. Let us recall that $v^{t}$ is such that, for every $i\in I$, 
\[
v^{t}\circ \mathrm{in}_{i}=v^{t}_{i}.
\tag{U4}\label{U4}
\]

Let $a$ be an element of $A$. Then the following chain of equalities holds
\allowdisplaybreaks
\begin{align*}
h(a)&=v^{t}(  h_{\overline{\Theta}}(a))
\tag{by Eq.~\ref{U3}}
\\&=
v^{t}(\mathrm{p}_{\Theta_{f(a)}}(a),f(a))
\tag{Def.~$h_{\overline{\Theta}}$, Prop.~\ref{PMult6}}
\\&=
v^{t}(\mathrm{in}_{f(a)}(\mathrm{p}_{\Theta_{f(a)}}(a)))
\tag{Def.~$\mathrm{in}_{f(a)}$}
\\&=
v^{t}_{f(a)}(\mathrm{p}_{\Theta_{f(a)}}(a))
\tag{by Eq.~\ref{U4}}
\\&=
\mathrm{in}_{t(f(a))}(v_{f(a)}(\mathrm{p}_{\Theta_{f(a)}}(a)))
\tag{Def.~$v^{t}_{f(a)}$, Prop.~\ref{PLaSumMor}}
\\&=
\mathrm{in}_{t(f(a))}(
v_{f(a)}(
\mathrm{pr}_{\Theta_{f(a)}}(
\mathrm{in}_{f(a),f(a)}(
a))))
\tag{Def.~$\mathrm{p}_{\Theta_{f(a)}}$, Prop.~\ref{PMult5}}
\\&=
\mathrm{in}_{t(f(a))}(
v_{f(a)}(
\mathrm{pr}_{\Theta_{f(a)}}(
a,f(a))))
\tag{Def.~$\mathrm{in}_{f(a),f(a)}$, Prop.~\ref{PMult2}}
\\&=
\mathrm{in}_{t(f(a))}(
v_{f(a)}(
[(a,f(a))]_{\Theta_{f(a)}}
))
\tag{Def.~$\mathrm{pr}_{\Theta_{f(a)}}$}
\\&=
(
v_{f(a)}(
[(a,f(a))]_{\Theta_{f(a)}}
),t(f(a))).
\tag{Def.~$\mathrm{in}_{t(f(a))}$, Prop.~\ref{PLaSumMor}}
\end{align*} 

From the last equality we conclude that, for every $a\in A$,
\begin{align*}
v_{f(a)}(
[(a,f(a))]_{\Theta_{f(a)}}
)&=\pi^{t(f(a))}_{0}(h(a))
.
\tag{U5}\label{U5}
\end{align*}

\begin{claim}\label{CMult710} 
The morphisms $v$ and $v_{\overline{\Theta}}$ are equal.
\end{claim}
We want to prove that, for every $j\in I$, $v_{j}=v_{\overline{\Theta},j}$.
Let $[(a,i)]_{\Theta_{j}}$ be an element of $M_{j}/{\Theta_{j}}$, where $i\in I$ is such that $(i,j)\in\leq_{\mathbf{I}}$ and $a\in A_{i}$. Then the following chain of equalities holds
\allowdisplaybreaks
\begin{align*}
v_{j}([(a,i)]_{\Theta_{j}})&=
v_{j}(
g^{\overline{\Theta}}_{i,j}(
[(a,i)]_{\Theta_{i}}
)
)
\tag{Prop.~\ref{PMult5}}
\\&=
g_{t(i),t(j)}(
v_{i}(
[(a,i)]_{\Theta_{i}}
)
)
\tag{By Eq.~\ref{U2}, $s=t$}
\\&=
g_{t(i),t(j)}(
\pi^{t(i)}_{0}(h(a))
)
\tag{By Eq.~\ref{U5}, $f$ is surj.}
\\&=
g_{t(i),t(j)}(
\pi^{t(i)}_{0}(h_{i}(a))
)
\tag{Def.~$h_{i}$}
\\&=
r_{i,j}(a)
\tag{By Eq.~\ref{E2}}
\\&=
r_{j}(\mathrm{in}_{i,j}(a))
\tag{By Eq.~\ref{E2}}
\\&=
r_{j}(a,i).
\tag{Def.~$\mathrm{in}_{i,j}$, Prop.~\ref{PMult2}}
\end{align*}

Thus, for every $j\in I$, the mapping $v_{j}$ is such that 
\begin{align*}
v_{j}\circ \mathrm{pr}_{\Theta_{j}}=r_{j}.
\label{U6}\tag{U6}
\end{align*}
But $v_{\overline{\Theta}, i}$ is the unique mapping from $M_{i}/{\Theta_{i}}$ to $N_{t(i)}$ such that 
$v_{\overline{\Theta}, i}\circ \mathrm{pr}_{\Theta_{i}}=r_{i}$. Therefore $v_{j} = v_{\overline{\Theta},j}$. 

This completes the proof of Claim~\ref{CMult710}.

All in all, we have that $s=s_{\overline{\Theta}}$ and $v=v_{\overline{\Theta}}$.

This completes the proof of Proposition~\ref{PMult7}.
\end{proof}

\begin{corollary}
The functor $\mathfrak{L}_{\Sigma}$ is a weak right multiadjoint.
\end{corollary}

\section{On the relationship between the P{\l}onka functor and the Lallement functor}

In this section, after recalling the essentials to sketch the definition of the 
P{\l}onka functor (for full details see~\cite{CC24}), we state the relationship between it and the Lallement functor.

To frame what follows, let us start by summarising the content of~\cite{CC24}. For a signature $\Sigma$ without $0$-ary operation symbols, two main categories are defined: the category $\mbox{\sffamily{\upshape{P{\l}Alg}}}(\Sigma)$, of P{\l}onka $\Sigma$-algebras, and the category $\int^{\mathsf{Sl}}\mathrm{Isys}_{\Sigma}$, of semilattice inductive systems of $\Sigma$-algebras. And the main result in~\cite{CC24} asserts that a certain functor, denoted by $\mathrm{Is}_{\Sigma}$, from 
$\mbox{\sffamily{\upshape{P{\l}Alg}}}(\Sigma)$ to $\int^{\mathsf{Sl}}\mathrm{Isys}_{\Sigma}$ has a left adjoint $\mbox{\upshape{P{\l}}}_{\Sigma}$, which is such that its object mapping assigns to a semilattice inductive systems of $\Sigma$-algebras its P{\l}onka sum.

\begin{assumption}
From now on $\Sigma$ stands for a plural signature, although we note that some of the following results are valid for signatures without $0$-ary operation symbols.
\end{assumption}

We begin by recalling the definition of the category of P{\l}onka $\Sigma$-algebras. 

\begin{definition}\label{DDOD4} 
Let $\mathbf{A}=(A,F)$ a $\Sigma$-algebra. A \emph{P{\l}onka operator for} $\mathbf{A}$ is a binary operation $D$ on $A$ such that, for every $x$, $y$, $z\in A$, every $n\in\mathbb{N}-\{0\}$, every $\sigma\in \Sigma_{n}$ and every $(x_{j})_{j\in n}\in A^{n}$, the following conditions are satisfied 
\begin{align*}
D(x,x)&=x;
\tag{D1}\label{D1}
\\
D(x,D(y,z))&=D(D(x,y),z);
\tag{D2}\label{D2}
\\
D(x,D(y,z))&=D(x,D(z,y)).
\tag{D3}\label{D3}
\\
D(F^{\mathbf{A}}_{\sigma}((x_{j})_{j\in n}), y)&=
F^{\mathbf{A}}_{\sigma}((D(x_{j},y))_{j\in n});
\tag{D4}\label{D4}
\\
D(y,F^{\mathbf{A}}_{\sigma}((x_{j})_{j\in n}))&=
D(y,D_{n}((x_{j})_{j\in n})).
\tag{D5}\label{D5}
\end{align*}
A \emph{P{\l}onka $\Sigma$-algebra} is an ordered pair $(\mathbf{A},D)$ in which $\mathbf{A}=(A,F)$ is a $\Sigma$-algebra and $D$ a P{\l}onka operator for 
$\mathbf{A}$ (see Remark~\ref{LnbD} below for the meaning of the operator $D_{n}$ in $\mathrm{D5}$). Given two P{\l}onka $\Sigma$-algebras $(\mathbf{A},D)$ and $(\mathbf{B},E)$, a \emph{morphism} $h$ from $(\mathbf{A},D)$ to $(\mathbf{B},E)$ is an ordered triple $((\mathbf{A},D),h,(\mathbf{B},E))$, abbreviated to $h\colon (\mathbf{A},D)\mor (\mathbf{B},E)$, where $h$ is a homomorphism from $\mathbf{A}$ to $\mathbf{B}$ and 
$ E\circ (h\times h)=h\circ D$. 
We denote by $\mbox{\sffamily{\upshape{P{\l}Alg}}}(\Sigma)$ the category whose objects are P{\l}onka $\Sigma$-algebras $(\mathbf{A},D)$ and whose morphisms are the morphisms $h\colon (\mathbf{A},D)\mor (\mathbf{B},E)$ of P{\l}onka $\Sigma$-algebras. 
\end{definition}

\begin{remark}
What we have called ``P{\l}onka operator for $\mathbf{A}$'' was called ``partition function for $\mathbf{A}$'' by  P{\l}onka in~\cite{P67}. Moreover, Definition~\ref{DDOD4}, taken from~\cite{RS85}, is the result of refining other definitions, such as Plonka's definitions in~\cite{{P67},{P74}}.
\end{remark}

\begin{remark}\label{LnbD}
In~\cite{CC24} we defined a \emph{left normal band} as an ordered pair $(A,D)$ in which $A$ is a set and $D$ a mapping $D\colon A\times A \mor A$ such that the above conditions $\mathrm{D1}$, $\mathrm{D2}$ and $\mathrm{D3}$ are satisfied. Moreover, for a left normal band $(A,D)$, we defined recursively both the family of the \emph{iterates of $D$ on the right}, $(D^{\mathrm{R}}_{n})_{n\geq 1}$, and the family of the \emph{iterates of $D$ on the left}, $(D^{\mathrm{L}}_{n})_{n\geq 1}$, and proved that, for every $n\in\mathbb{N}-\{0\}$ and $(x_{j})_{j\in n}\in A^{n}$,
$
D^{\mathrm{R}}_{n}((x_{j})_{j\in n})=D^{\mathrm{L}}_{n}((x_{j})_{j\in n}).
$ 
This making it legitimate to use $D_{n}$ instead of $D^{\mathrm{R}}_{n}$ or $D^{\mathrm{L}}_{n}$.
\end{remark}


We next state that for a P{\l}onka $\Sigma$-algebra  $(\mathbf{A},D)$, the P{\l}onka operator $D$ induces a congruence on $(\mathbf{A},D)$.

\begin{proposition}\label{PDOD4Cong} 
Let $(\mathbf{A},D)$ be a P{\l}onka $\Sigma$-algebra and $\Phi^{D}$ the binary relation on $A$ defined, for every $x,y\in A$, as: $(x,y)\in \Phi^{D}$ if, and only if, $D(x,y)=x$ and $D(y,x)=y$. Then $\Phi^{D}$ is a congruence on $(\mathbf{A},D)$, i.e., a congruence on $\mathbf{A}$ which, in addition, is compatible with $D$. We will call $\Phi^{D}$ the \emph{congruence induced by $D$ on $(\mathbf{A},D)$}.
\end{proposition}

In the following, we establish the relationship between the Lallement systems of 
$\Sigma$-algebras labeled by an npi $\Sigma$-algebra and the inductive systems of 
$\Sigma$-algebras labeled by a semilattice. To do it, we start by establishing the relationship between the category $\mathsf{Sl}$ and the category $\mathsf{NPIAlg}(\Sigma)$ as well as between some categories derived from them.

\begin{proposition}\label{SlNPIAlg}
There exists a functor $G$ from $\mathsf{Sl}$ to $\mathsf{NPIAlg}(\Sigma)$ and a natural transformation $\eta$ from $\mathrm{Isys}_{\Sigma}$ to $\mathrm{Isys}_{\Sigma}\circ G^{\mathrm{op}}$ as well as a functor, also denoted by $G$, from the category $\mathsf{Sl}_{\mathsf{s}}$, of semilattices and surjective homomorphisms between semilattices, to the category 
$\mathsf{NPIAlg}(\Sigma)_{\mathsf{s}}$ and a natural transformation, also denoted by $\eta$, from $\mathrm{Isys}_{\Sigma}$ to $\mathrm{Isys}_{\Sigma}\circ G^{\mathrm{op}}$. 
\end{proposition}

\begin{proof}
Given a semilattice $\mathbf{I}$, it suffices to define $G(\mathbf{I})$ as the npi 
$\Sigma$-algebra whose underlying set is $I$ and whose structural operations are defined as follows: for every $n\in \mathbb{N}-\{0\}$ and every $\sigma\in \Sigma_{n}$, the $n$-ary structural operation $F_{\sigma}^{\mathbf{I}}$ on $I$ associated to $\sigma$ sends  $(i_{k})_{k\in n}\in I^{n}$ to $F_{\sigma}^{\mathbf{I}}((i_{k})_{k\in n}) = \bigvee_{k\in n}i_{k}$. The morphism mapping of $G$ assigns to a morphism $h$ from a semilattice $\mathbf{I}$ to another $\mathbf{P}$ the same morphism, but considered as a homomorphism from $G(\mathbf{I})$ to $G(\mathbf{P})$. This last  definition is sound because, for every $n\in \mathbb{N}-\{0\}$, every $\sigma\in \Sigma_{n}$ and every $(i_{k})_{k\in n}\in I^{n}$, $h(\bigvee_{k\in n} i_{k}) = \bigvee_{k\in n}h(i_{k})$.

With regard to the natural transformation $\eta$ from $\mathrm{Isys}_{\Sigma}$ to 
$\mathrm{Isys}_{\Sigma}\circ G^{\mathrm{op}}$, suffice it to say that 
every $\mathbf{I}$-inductive system of $\Sigma$-algebras  
$
\mathcal{A} = ((\mathbf{A}_{i})_{i\in I},(f_{i,j})_{(i,j)\in \leq})
$  
has canonically associated with it an $\mathbf{I}$-Lallement system of $\Sigma$-algebras:  
$\mathcal{A}=((\mathbf{A}_{i}, \mathbf{A}_{i})_{i\in I}, (f_{i,j})_{(i,j)\in\leq_{\mathbf{I}}})$ (and the same applies to the morphisms).
\end{proof}

\begin{corollary}
There exists a functor $\int^{G}\eta$ from 
$\int^{\mathsf{Sl}}\mathrm{Isys_{\Sigma}}$ to $\int^{\mathsf{NPIAlg}(\Sigma)}\mathrm{Lsys}_{\Sigma}$ such that $G\circ \pi_{\mathrm{Isys}_{\Sigma}} = \pi_{\mathrm{Lsys}_{\Sigma}}\circ \int^{G}\eta$ as well as a functor, which we shall continue to denote by $\int^{G}\eta$, from $\int^{\mathsf{Sl}_{\mathsf{s}}}\mathrm{Isys_{\Sigma}}$ to $\int^{\mathsf{NPIAlg}(\Sigma)_{\mathsf{s}}}\mathrm{Lsys}_{\Sigma}$ with the same property.
\end{corollary}

\begin{proof}
By Proposition~\ref{SlNPIAlg} and Proposition~\ref{Intmorf}. 
\end{proof}

We continue by defining the notions of semilattice inductive system of $\Sigma$-algebras relative to a semilattice and of morphism between them. 

\begin{definition}\label{DIndSys}
Let $\mathbf{I}$ be a semilattice. Then a \emph{semilattice inductive system of 
$\Sigma$-algebras relative to} $\mathbf{I}$ or an $\mathbf{I}$-\emph{inductive system of $\Sigma$-algebras} for brevity, is an ordered pair 
$
\mathcal{A} = ((\mathbf{A}_{i})_{i\in I},(f_{i,j})_{(i,j)\in \leq}),
$ 
where $(\mathbf{A}_{i})_{i\in I}$ is an $I$-indexed family of $\Sigma$-algebras, and 
$(f_{i,j})_{(i,j)\in \leq}$ a family of homomorphisms in $\prod_{(i,j)\in \leq}\mathrm{Hom}(\mathbf{A}_{i},\mathbf{A}_{j})$ such that, for every $i\in I$, $f_{i,i} = \mathrm{id}_{\mathbf{A}_{i}}$ and, for every 
$i$, $j$, $k\in I$, if $i\leq j$ and $j\leq k$, then $f_{j,k}\circ f_{i,j} = f_{i,k}$, i.e., a covariant functor from the category canonically associated to $\mathbf{I}$ to 
$\mathsf{Alg}(\Sigma)$. Let $\mathcal{A}$ and $\mathcal{B} = ((\mathbf{B}_{i})_{i\in I},(g_{i,j})_{(i,j)\in \leq})$ be $\mathbf{I}$-inductive systems of 
$\Sigma$-algebras. Then a \emph{morphism} from $\mathcal{A}$ to $\mathcal{B}$ is an ordered triple $(\mathcal{A},u,\mathcal{B})$, abbreviated to $u\colon \mathcal{A}\mor \mathcal{B}$, where $u = (u_{i})_{i\in I}$ is a family of homomorphisms in $\prod_{i\in I}\mathrm{Hom}(\mathbf{A}_{i},\mathbf{B}_{i})$ such that, for every $(i,j)\in \leq$, $u_{j}\circ f_{i,j} = g_{i,j}\circ u_{i}$, i.e., a natural transformation from $\mathcal{A}$ to $\mathcal{B}$. We denote by $\mathsf{Alg}(\Sigma)^{\mathbf{I}}$ the corresponding category. 
\end{definition}

After this we define the assignment $\mathrm{Isys_{\Sigma}}$ from $\mathsf{Sl}$ to $\mathsf{Cat}$.

\begin{definition}\label{DIsys}
We let $\mathrm{Isys_{\Sigma}}$ stand for the assignment from $\mathsf{Sl}$ to $\mathsf{Cat}$ defined as follows:
\begin{enumerate}
\item for every semilattice $\mathbf{I}$, $\mathrm{Isys_{\Sigma}}(\mathbf{I})$ is $\mathsf{Alg}(\Sigma)^{\mathbf{I}}$, and 
\item for every morphism $\xi\colon\mathbf{I}\mor \mathbf{P}$, $\mathrm{Isys_{\Sigma}}(\xi)$ is the functor  from $\mathsf{Alg}(\Sigma)^{\mathbf{P}}$ to $\mathsf{Alg}(\Sigma)^{\mathbf{I}}$ defined as follows:
\begin{enumerate}
\item for every $\mathbf{P}$-inductive systems of $\Sigma$-algebras $\mathcal{A} = ((\mathbf{A}_{p})_{p\in P},(f_{p,q})_{(p,q)\in \leq})$, $\mathrm{Isys}_{\Sigma}(\xi)(\mathcal{A})$, denoted by $\mathcal{A}_{\xi}$ for short,  is $((\mathbf{A}_{\xi(i)})_{i\in I},(f_{\xi(i),\xi(j)})_{(i,j)\in \leq})$, which is an $\mathbf{I}$-inductive systems of $\Sigma$-algebras, and
\item for every morphism of $\mathbf{P}$-inductive systems of $\Sigma$-algebras $u = (u_{p})_{p\in P}$ from $\mathcal{A}$ to $\mathcal{B} = ((\mathbf{B}_{p})_{p\in P},(g_{p,q})_{(p,q)\in \leq})$, 
$\mathrm{Isys}_{\Sigma}(\xi)(u)$, denoted by $u_{\xi}$ for short, is $(u_{\xi(i)})_{i\in I}$, which is a morphism of $\mathbf{I}$-inductive systems of $\Sigma$-algebras from $\mathcal{A}_{\xi}$ to $\mathcal{B}_{\xi}$.
\end{enumerate} 
\end{enumerate}
We also let $\mathrm{Isys_{\Sigma}}$ stand for the restriction of 
$\mathrm{Isys_{\Sigma}}$ to $\mathsf{Sl}_{\mathsf{s}}$.
\end{definition}


\begin{proposition}
The assignment $\mathrm{Isys_{\Sigma}}$ from $\mathsf{Sl}$ to $\mathsf{Cat}$ as well as its restriction to $\mathsf{Sl}_{\mathsf{s}}$ is a contravariant functor.
\end{proposition}

\begin{definition}\label{DIntSsl}
From the split indexed category $(\mathsf{Sl},\mathrm{Isys_{\Sigma}})$ we obtain, by means of the Grothendieck construction, the category $\int^{\mathsf{Sl}}\mathrm{Isys_{\Sigma}}$ which has as objects the \emph{semilattice inductive system of $\Sigma$-algebras} or, simply, the  \emph{inductive system of $\Sigma$-algebras}, i.e., the ordered pairs 
$
\mathbcal{A} = (\mathbf{I},\mathcal{A}),
$
where $\mathbf{I}$ is a semilattice and $\mathcal{A}$ an $\mathbf{I}$-inductive system of 
$\Sigma$-algebras; and as morphisms from $\mathbcal{A}$ to $\mathbcal{B} = (\mathbf{P},\mathcal{B})$ the ordered triples $(\mathbcal{A},(\xi,u),\mathbcal{B})$, abbreviated to $(\xi,u)\colon \mathbcal{A}\mor \mathbcal{B}$, where $\xi$ is a morphism from $\mathbf{I}$ to $\mathbf{P}$ and $u$ a morphism from the $\mathbf{I}$-inductive system of $\Sigma$-algebras $\mathcal{A}$ to the $\mathbf{I}$-inductive system of $\Sigma$-algebras $\mathcal{B}_{\xi}$. Moreover, we let 
$\pi_{\mathrm{Isys_{\Sigma}}}$ stand for the canonical split fibration, i.e., the functor from $\int^{\mathsf{Sl}}\mathrm{Isys_{\Sigma}}$ to $\mathsf{Sl}$ that sends $\mathbcal{A} = (\mathbf{I},\mathcal{A})$ to $\mathbf{I}$ and $(\xi,u)\colon \mathbcal{A}\mor \mathbcal{B}$, with $\mathbcal{B} = (\mathbf{P},\mathcal{B})$, to $\xi\colon \mathbf{I}\mor \mathbf{P}$. In the same way, from the split indexed category $(\mathsf{Sl}_{\mathsf{s}},\mathrm{Isys_{\Sigma}})$ we obtain the category 
$\int^{\mathsf{Sl}_{\mathsf{s}}}\mathrm{Isys_{\Sigma}}$ and the canonical split fibration, also denoted by $\pi_{\mathrm{Isys_{\Sigma}}}$, from $\int^{\mathsf{Sl}_{\mathsf{s}}}\mathrm{Isys_{\Sigma}}$ to $\mathsf{Sl}_{\mathsf{s}}$.
\end{definition}

\begin{proposition}
There are functors $\mathrm{Is}_{\Sigma}$ from 
$\mbox{\sffamily{\upshape{P{\l}Alg}}}(\Sigma)$ to $\int^{\mathsf{Sl}}\mathrm{Isys}_{\Sigma}$ and $\mbox{\upshape{P{\l}}}_{\Sigma}$, the \emph{P{\l}onka functor}, from 
$\int^{\mathsf{Sl}}\mathrm{Isys}_{\Sigma}$ to $\mbox{\sffamily{\upshape{P{\l}Alg}}}(\Sigma)$ such that $\mbox{\upshape{P{\l}}}_{\Sigma}$ is a left adjoint to $\mathrm{Is}_{\Sigma}$. We also let $\mbox{\upshape{P{\l}}}_{\Sigma}$ stand for the birestriction of $\mbox{\upshape{P{\l}}}_{\Sigma}$ to $\int^{\mathsf{Sl}_{\mathsf{s}}}\mathrm{Isys}_{\Sigma}$ and $\mbox{\sffamily{\upshape{P{\l}Alg}}}(\Sigma)_{\mathsf{s}}$.
\end{proposition}

\begin{proof}
See~\cite{CC24} Propositions~6.14 and~6.16. Suffice it to say in this regard and for our objetive in this section that the object mapping of $\mbox{\upshape{P{\l}}}_{\Sigma}$ assigns to an inductive system of $\Sigma$-algebras $\mathbcal{A} = (\mathbf{I},\mathcal{A})$, where the $\mathbf{I}$-inductive system of $\Sigma$-algebras $\mathcal{A}$ is $((\mathbf{A}_{i})_{i\in I},(f_{i,j})_{(i,j)\in \leq})$, the  P{\l}onka $\Sigma$-algebra 
$(\mbox{\bfseries{\upshape{P{\l}}}}_{\Sigma}(\mathbcal{A}), D^{\mathbcal{A}})$, where  the $\Sigma$-algebra 
$\mbox{\bfseries{\upshape{P{\l}}}}_{\Sigma}(\mathbcal{A})$, the \emph{P{\l}onka sum} of $\mathbcal{A}$, is defined as follows:
\begin{enumerate}
\item the underlying set of $\mbox{\bfseries{\upshape{P{\l}}}}_{\Sigma}(\mathbcal{A})$ is 
\[
\textstyle
\coprod_{i\in I}A_{i} = \bigcup_{i\in I}(A_{i}\times \{i\}),\] 
\item for every $n\in\mathbb{N}-\{0\}$ and every $\sigma\in \Sigma_{n}$, the structural operation $F^{\scriptsize{\mbox{\bfseries{\upshape{P{\l}}}}_{\Sigma}(\mathbcal{A})}}_{\sigma}$ is defined as
\[
\textstyle
F^{\scriptsize{\mbox{\bfseries{\upshape{P{\l}}}}_{\Sigma}(\mathbcal{A})}}_{\sigma}
\left\lbrace
\begin{array}{ccc}
(\bigcup_{i\in I}(A_{i}\times\{i\}))^{n}
&\mor&
\bigcup_{i\in I} (A_{i}\times\{i\})
\\
(x_{j},i_{j})_{j\in n}
&\longmapsto&
(
F^{\mathbf{A}_{\bigvee_{j\in n}i_{j}}}_{\sigma}(
(f_{i_{j},\bigvee_{j\in n}i_{j}}(
x_{j}
))_{j\in n}),
\bigvee_{j\in n}i_{j}
)
\end{array}
\right.
\]
\end{enumerate} 
and the P{\l}onka operator $D^{\mathbcal{A}}$ is the mapping defined as follows:
\[
D^{\mathbcal{A}}
\left\lbrace
\begin{array}{ccc}
 \bigcup_{i\in I}\left(A_{i}\times\{i\}\right)\times \bigcup_{i\in I}\left(A_{i}\times\{i\}\right)&\mor&\bigcup_{i\in I}\left(A_{i}\times\{i\}\right)
 \\
 ((x,j),(y,k))&\longmapsto&
 \left(
f_{j,j\vee k}(x),j\vee k
\right)
\end{array}
\right.
\]
\end{proof}

\begin{definition}
We let $P$ stand for the functor from $\mbox{\sffamily{\upshape{P{\l}Alg}}}(\Sigma)_{\mathsf{s}}$ to $\mathsf{Alg}(\Sigma)_{\mathsf{s}}$ that sends 
$(\mathbf{A},D)$ to $\mathbf{A}/\Phi^{D}$; and $Q$ for the functor from $\int_{\mathsf{NPIAlg}(\Sigma)_{\mathsf{s}}}(\mathsf{Alg}(\Sigma)\,{\downarrow_{\mathsf{s}}}\, \bigcdot)$ to  
$\mathsf{Alg}(\Sigma)_{\mathsf{s}}$ that sends $(\mathbf{I},(\mathbf{A},f))$ to $\mathbf{I}$.
\end{definition} 

After these preliminaries, finally, we state the relationship between the P{\l}onka functor and the Lallement functor.

\begin{proposition}
There exists a natural isomorphism $\theta$ from the functor $P\circ \mbox{\upshape{P{\l}}}_{\Sigma}$ to the functor $Q\circ\mathfrak{L}_{\Sigma}\circ\int^{G}\eta$ (this is illustrated in \emph{Figure~\ref{FigPl}}). 
\end{proposition}

\begin{proof}
For every $\mathbcal{A}$ in $\int^{\mathsf{Sl}_{\mathsf{s}}}\mathrm{Isys_{\Sigma}}$, the underlying Lallement sum of $\mathfrak{L}_{\Sigma}(\int^{G}\eta(\mathbcal{A}))$ coincides with the underlying P{\l}onka sum of $\mbox{\upshape{P{\l}}}_{\Sigma}(\mathbcal{A})$ and $\mathrm{p}_{\mathbcal{A}}$ is surjective. Moreover, since
$\mathrm{Ker}(\mathrm{p}_{\mathbcal{A}}) = \Phi^{D^{\mathbcal{A}}}$, because, for every $(x,j)$, $(y,k)\in \mbox{\upshape{P{\l}}}_{\Sigma}(\mathbcal{A})$, $((x,j),(y,k))\in \Phi^{D^{\mathbcal{A}}}$ if, and only if, $j = k$ (by Proposition 6.16, and Claim 6.22 in~\cite{CC24}) and $((x,j),(y,k))\in \mathrm{Ker}(\mathrm{p}_{\mathbcal{A}})$ if, and only if, $j = k$ (by definition of $\mathrm{p}_{\mathbcal{A}}$), we have that there exists a unique natural isomorphism $\theta_{\mathbcal{A}}$ from 
$\mbox{\bfseries{\upshape{P{\l}}}}_{\Sigma}(\mathbcal{A})/\Phi^{D^{\mathbcal{A}}}$ to  $\mathbf{I}$ such that 
$\theta_{\mathbcal{A}}\circ\mathrm{pr}_{\Phi^{D^{\mathbcal{A}}}} = \mathrm{p}_{\mathbcal{A}}$. Thus it suffices to define the mapping 
$\theta$ from the set of objects of $\int^{\mathsf{Sl}_{\mathsf{s}}}\mathrm{Isys_{\Sigma}}$ to the set of morphisms of $\mathsf{Alg}(\Sigma)_{\mathsf{s}}$ as the family $(\theta_{\mathbcal{A}})_{\mathbcal{A}\in\mathrm{Ob}(\int^{\mathsf{Sl}_{\mathsf{s}}}\mathrm{Isys_{\Sigma}})}$.
\end{proof}

\begin{figure}
\centering
\begin{tikzpicture}
[AClimentM/.style={|-{To [angle'=45, length=5.75pt, width=4pt, round]}},
ACliment/.style={-{To [angle'=45, length=5.75pt, width=4pt, round]}},scale=1,
AClimentD/.style={double equal sign distance,-implies}, scale=1]
            
\node[] (P) 		at (0,0) [] {$\mbox{\sffamily{\upshape{P{\l}Alg}}}(\Sigma)_{\mathsf{s}}$};
\node[] (I1) 	at (0,-2) [] {$\int^{\mathsf{Sl}_{\mathsf{s}}}\mathrm{Isys}_{\Sigma}$};
\node[] (I2) 	at (7,-2) [] {$\int^{\mathsf{NPIAlg}(\Sigma)_{\mathsf{s}}}\mathrm{Lsys}_{\Sigma}$};
\node[] (A) 		at (3,0) [] {$\mathsf{Alg}(\Sigma)_{\mathsf{s}}$};
\node[] (I3) 	at (7,0) [] {$\int_{\mathsf{NPIAlg}(\Sigma)_{\mathsf{s}}}(\mathsf{Alg}(\Sigma)\,{\downarrow_{\mathsf{s}}}\, \bigcdot)$};
            
\draw[ACliment]  (P) to node [above]	{$P$} (A);
\draw[ACliment]  (I1) to node [left]	{$\mbox{\upshape{P{\l}}}_{\Sigma}$} (P);
\draw[ACliment]  (I1) to node [below, pos=0.42]	{$\int^{G}\eta$} (I2);
\draw[ACliment]  (I2) to node [right]	{$\mathfrak{L}_{\Sigma}$} (I3);
\draw[ACliment]  (I3) to node [above]	{$Q$} (A);

\node[] (auxsc) at (2.25,-.65) [] {};
\node[] (auxtg) at (3.75,-1.15) [] {};

\draw[AClimentD]  (auxsc) to node [above]	{$\theta$} (auxtg);
            
\end{tikzpicture}
\caption{The natural isomorphism $\theta$.}
\label{FigPl}
\end{figure}

\section*{Acknowledgements}
The second author was supported by the grant CIAICO/2023/007 from the
Conselleria d'Educaci\'{o}, Universitats i Ocupaci\'{o}, Generalitat Valenciana.
The authors sincerely thank the reviewer for their valuable suggestions and the editor for their prompt handling of the manuscript.

\end{document}